\documentclass[a4paper]{amsart}

\usepackage{amssymb}
\usepackage{amsthm}
\usepackage{latexsym}

\usepackage{color}%
\definecolor{darkmagenta}{rgb}{.5, 0, .5}

\newtheorem{thm}{Theorem}[section]
\newtheorem{prop}[thm]{Proposition}
\newtheorem{lemma}[thm]{Lemma}
\newtheorem{cor}[thm]{Corollary}
\theoremstyle{definition}

\newtheorem{ex}[thm]{Example}
\newtheorem{rem}[thm]{Remark}

 \DeclareMathOperator{\Spec}{Spec}
 \DeclareMathOperator{\Max}{Max}
 \DeclareMathOperator{\Inv}{Inv}
 \DeclareMathOperator{\Na}{Na}
 
 \DeclareMathOperator{\QInv}{QInv}
\DeclareMathOperator{\QSpec}{QSpec}
 \DeclareMathOperator{\QMax}{QMax}

\begin{document}

\title[SEMISTAR INVERTIBILITY]{SEMISTAR INVERTIBILITY ON INTEGRAL DOMAINS}
\author{Marco Fontana \hskip 1.8cm Giampaolo Picozza}



\date{\today}

\dedicatory{\textrm{Dipartimento di Matematica, Universit\`a degli Studi
``Roma Tre'', \\
\vskip 0.1cm Largo San Leonardo Murialdo 1, 00146 Roma }\\
\rm
\vskip 0.2cm
\texttt{fontana@mat.uniroma3.it \;\;\; \;\;\;\;  picozza@mat.uniroma3.it}}
%
\begin{abstract} After the introduction in 1994, by Okabe and Matsuda, of the
notion of semistar operation, many authors have investigated different
aspects of this general and powerful concept.  A natural development of
the recent work in this area leads to investigate the concept of
invertibi\-li\-ty in the semistar setting.  In this paper, we will show
the existence of a ``theoretical obstruction'' for extending many
results, proved for star-invertibility, to the semistar case.  For this
reason, we will introduce two distinct notions of invertibility in the
semistar setting (called $\star$--invertibility and
quasi--$\star$--invertibility), we will discuss the motivations of these
``two levels'' of invertibility and we will extend, accordingly, many
classical results proved for the $d$--, $v$--, $t$-- and $w$--\
invertibility.  Among the main properties proved here, we mention the
following:\ \bf (a) \rm several characterizations of
$\star$--invertibility and quasi--$\star$--invertibility and necessary
and sufficient conditions for the equivalence of these two notions; \
\bf (b) \rm the relations between the $\star$--invertibility (or
quasi--$\star$--invertibility) and the invertibility (or
quasi--invertibility) with respect to the semistar operation of finite
type (denoted by $\star_{_{\!f}}$) and to the stable semistar operation
of finite type (denoted by $\widetilde{\star}$), canonically associated
to $\star$;\ \bf (c) \rm a characterization of the H$(\star)$--domains
in terms of semistar--invertibility (note that the H$(\star)$--domains
generalize, in the semistar setting, the H--domains introduced by Glaz
and Vasconcelos);\ \bf (d) \rm for a semistar operation of finite type a
nonzero finitely generated (fractional) ideal $I$ is $\star$--invertible
(or, equivalently, quasi--$\star$--invertible, in the stable semistar
case) if and only if its extension to the Nagata semistar ring
$I\Na(D,\star)$ is an invertible ideal of $\Na(D,\star)$.
\end{abstract}
\maketitle

\section{ Introduction and background results}

  The notions of $t$--invertibility, $v$--invertibility and $w$--invertibility, that generalize the
classical concept of ($d$--)invertibility (these definitions will be
recalled in Section 2), have been introduced in the recent years for a
better understanding of the multiplicative (ideal) structure of integral
domains.  In particular, $t$--invertibility has a key role for extending
the notion of class group from Krull domains to general integral domains
(cf.  \cite{B82}, \cite{B83}, \cite {BZ88} and the survey paper
\cite{A00}).  An interesting chart of a large set of various $t$--,
$v$--, $d$-- invertibility based characterizations of several classes of
integral domains can be found at the end of \cite{AMZ89}; some
motivations for introducing the $w$--invertibility and the first
properties showing the ``good'' behaviour of this notion can be found in
\cite{WM97} (cf.  also \cite{HH80}).  The concept of star operation (or,
equivalently, ideal system, cf.  the books by Jaffard
\cite{Jaffard60}, Gilmer \cite{Gil92} and Halter-Koch \cite{HK98})
provides an abstract setting for approaching these different aspects of
invertibility.  A recent paper by Zafrullah \cite{Z2} gives an excellent
and updated survey of this point of view.

After the introduction in 1994, by Okabe and Matsuda \cite{OM94}, of the
notion of semistar operation, as a more general and natural setting for
studying multiplicative systems of ideals and modules, many authors have
investigated the possible extensions to the semistar setting of different
aspects of the classical theory of ideal systems, based on the
pioneering work by W. Krull, E. Noether, H. Pr\"ufer and P. Lorenzen from
1930's (cf.  for instance \cite {MSu95}, \cite{MS96}, \cite{OM97},
\cite{M98}, \cite{M00}, \cite{FH2000}, \cite{HK01}, \cite{FL01a},
\cite{FL01b}, \cite{FL03}, \cite{ElBF03}, \cite{FJS03}, \cite{FJS03b},
\cite{HK03}, \cite{Mi03}, \cite{MiSa03}, \cite{O03}).

A natural development of this work leads to investigate the concept of
invertibi\-li\-ty in the semistar setting.  This is the purpose of the
present paper, in which we will show the existence of a ``theoretical
obstruction'' for extending many results, proved for star-invertibility,
to the semistar case.  For this reason, we will be forced to introduce two
distinct notions of invertibility in the semistar setting (called
$\star$--invertibility and quasi--$\star$--invertibility; the explicit
definitions will be given in Section 2), we will discuss the motivations
of these ``two levels'' of invertibility and we will extend,
accordingly, many classical results proved for the  $d$--, $v$--,
$t$-- and $w$--\ invertibility.

Among the main properties proved in this work, we mention
the following:\ \bf (a) \rm several characterizations of
$\star$--invertibility and quasi--$\star$--invertibility and necessary
and sufficient conditions for the equivalence of these two notions; \
\bf (b) \rm the relations between the $\star$--invertibility (or
quasi--$\star$--invertibility) and the invertibility (or
quasi--invertibility) with respect to the semistar operation of finite
type (denoted by $\star_{_{\!f}}$) and to the stable semistar operation
of finite type (denoted by $\widetilde{\star}$), canonically associated
to $\star$\ [in case, $\star =v$ is the Artin's $v$--operation, then
$\star_{_{\!f}} =t$ and $\widetilde{\star}= w$];\ \bf (c) \rm a
characterization of the H$(\star)$--domains in terms of
semistar--invertibility (note that the H$(\star)$--domains generalize in
the semistar setting the H--domains introduced by Glaz and Vasconcelos
\cite{GV77}, more precisely, we will see in Section 2 that an H--domain
coincides with an H$(v)$--domain);\ \bf (d) \rm for a semistar operation
of finite type a nonzero finitely generated (fractional) ideal $I$ is
$\star$--invertible (or, equivalently, quasi--$\star$--invertible, in
the stable semistar case) if and only if its extension to the Nagata
semistar ring $I\Na(D,\star)$ is an invertible ideal of $\Na(D,\star)$
(the definition of $\Na(D,\star)$ will be recalled at the end of this
section).

\medskip
\centerline {$\boldsymbol{\ast \, \ast \, \ast\, \ast \, \ast }$}

\medskip

Let $D$ be an integral domain with quotient field $K$. Let
$\boldsymbol{\overline{F}}(D)$ denote the set of all nonzero
$D$--submodules of $K$ and let $\boldsymbol{F}(D)$ be the set of
all nonzero fractional ideals of $D$, i.e. $E \in
\boldsymbol{F}(D)$ if $E \in \boldsymbol{\overline{F}}(D)$ and
there exists a nonzero $d \in D$ with $dE \subseteq D$. Let
$\boldsymbol{f}(D)$ be the set of all nonzero finitely generated
$D$--submodules of $K$. Then, obviously $\boldsymbol{f}(D)
\subseteq \boldsymbol{F}(D) \subseteq
\boldsymbol{\overline{F}}(D)$.

A map $\star: \boldsymbol{\overline{F}}(D) \to
\boldsymbol{\overline{F}}(D), E \mapsto E^\star$, is called a
\emph{semistar operation} on $D$ if, for all $x \in K$, $x \neq 0$, and
for all $E,F \in \boldsymbol{\overline{F}}(D)$, the following properties
hold:
\begin{enumerate}
\item[$(\star_1)$] $(xE)^\star=xE^\star$;
\item[$(\star_2)$] $E \subseteq F$ implies $E^\star \subseteq
F^\star$;
\item[$(\star_3)$] $E \subseteq E^\star$ and $E^{\star \star} :=
\left(E^\star \right)^\star=E^\star$.
\end{enumerate}
cf.  for instance \cite{FH2000}.  Recall that \cite[Theorem 1.2 and p.
174]{FH2000}, for all $E,F \in \boldsymbol{\overline{F}}(D)$,
we have : $$\begin{array}{rl}
(EF)^\star =& \hskip -7pt (E^\star F)^\star =\left(EF^\star\right)^\star
=\left(E^\star F^\star\right)^\star\,;\\
(E+F)^\star =& \hskip -7pt \left(E^\star + F\right)^\star= \left(E +
F^\star\right)^\star= \left(E^\star + F^\star\right)^\star\,;\\
(E:F)^\star \subseteq & \hskip -7pt (E^\star :F^\star) = (E^\star :F) =
\left(E^\star :F\right)^\star\,;\\
(E\cap F)^\star \subseteq & \hskip -7pt E^\star \cap F^\star = \left(E^\star
\cap F^\star \right)^\star\,,\;\, \mbox{ if $E\cap F \neq (0)$}\,;
\end{array}
$$

When $D^\star=D$, we say that $\star$ is a \emph{(semi)star
operation}, since, restricted to $\boldsymbol{F}(D)$ it is a star
operation.

For star operations it is very useful the notion of
\emph{$\star$--ideal}, that is, a  nonzero ideal $I \subseteq D$, such that
$I^\star = I$.  For semistar operations we need a more general notion,
that coincides with the notion of $\star$--ideal, when $\star$ is a
(semi)star operation.  We say that a  nonzero (integral) ideal $I$ of
$D$ is a \emph{quasi--$\star$--ideal} if $I^\star \cap D = I$.  For
example, it is easy to see that, if $I^\star \neq D^\star$, then
$I^\star \cap D$ is a quasi--$\star$--ideal that contains $I$ (in
particular, a $\star$--ideal is a quasi--$\star$--ideal).  Note that
$I^\star \neq D^\star$ is equivalent to $I^\star \cap D \neq D$.

A \emph{quasi--$\star$--prime} is a quasi--$\star$--ideal which is also a
prime ideal.  We call a \emph{quasi--$\star$--maximal} a maximal element
in the set of all proper quasi--$\star$--ideals of $D$.  We denote by
$\QSpec^\star(D)$ (respectively, $\QMax^\star(D)$) the set of all
quasi--$\star$--prime (respectively, quasi--$\star$--maximal).

If $\star$ is a semistar operation on $D$, then we can consider a
map\ $\star_{\!_f}:  \boldsymbol{\overline{F}}(D) \to
\boldsymbol{\overline{F}}(D)$ defined for each $E \in
\boldsymbol{\overline{F}}(D)$ as follows: $E^{\star_{\!_f}}:=\bigcup
\{F^\star \, \vert \, F \in \boldsymbol{f}(D) \mbox{ and } F
\subseteq E\}$.  It is easy to see that $\star_{\!_f}$ is a semistar
operation on $D$, called \emph{the semistar operation of finite type
associated to $\star$}.  Note that, for each $F \in \boldsymbol{f}(D)$,
$F^\star=F^{\star_{\!_f}}$.  A semistar operation $\star$ is called a
\emph{semistar operation of finite type} if $\star=\star_{\!_f}$.  It is easy
to see that $(\star_{\!_f}\!)_{\!_f}=\star_{\!_f}$ (that is, $\star_{\!_f}$ is of finite
type).

 If ${\star_1}$ and ${\star_2}$ are
two semistar operations on $D$, we say that ${\star_1} \leq
{\star_2}$ if $E^{\star_1} \subseteq E^{\star_2}$, for each $E \in
\boldsymbol{\overline{F}}(D)$.   In this situation, it is easy to
see that $\left(E^{\star_{1}}\right)^{\star_{2}} = E^{\star_2}=
\left(E^{\star_{2}}\right)^{\star_{1}}$.  Obviously, for each semistar
operation $\star$, we have $\star_{\!_f} \leq \star$.

The following result, with a different terminology, was proved in
\cite{FH2000} (cf. also \cite[Lemma 2.3]{FL03}).

\begin{lemma}
Let $\star$ be a semistar operation on an integral domain $D$.
Assume that $\star$ is non trivial and $\star=\star_{\!_f}$. Then:
\begin{enumerate}
\item[(1)] Each proper quasi--$\star$--ideal is contained in a
quasi--$\star$--maximal.  \item[(2)] Each quasi--$\star$--maximal is a
quasi--$\star$--prime.  \item[(3)] Set $\Pi^\star := \{P \in \Spec(D) \,
\vert \, P \neq 0 \mbox{ and } P^\star \cap D \neq D \}$, then
$\QSpec^\star(D) \subseteq \Pi^\star$ and the set of maximal elements
$\Pi^{\star}_{\mbox{\rm \tiny \textsl{max}}}$ of $\Pi^\star$ is nonempty
and coincides with $\QMax^\star(D)$.  \hfill $\Box$
\end{enumerate}
\end{lemma}

For the sake of simplicity,  we will denote simply by
$\mathcal{M}(\star)$ the set $\QMax^\star(D)$ of the
quasi--$\star$--maximal ideals of $D$.

If $\Delta \subseteq \Spec(D)$, the map $\star_{\Delta}:
\boldsymbol{\overline{F}}(D) \mapsto
\boldsymbol{\overline{F}}(D)$, $E  \mapsto E^{\star_\Delta} := \bigcap
\{ED_P \, \vert \, P \in \Delta \}$, is a semistar operation.  If
$\star=\star_\Delta$, for some $\Delta \subseteq \Spec(D)$, we say that
$\star$ is a \emph{spectral semistar operation}.  In particular, if
$\Delta=\{P\}$, then $\star_{\{P\}}$ is the semistar operation on $D$
defined by $E^{\star_{\{P\}}}:=ED_P$, for each $E \in
\boldsymbol{\overline{F}}(D)$.  We say that a \emph{semistar operation}
is \emph{stable} if $(E \cap F)^\star= E^\star \cap F^\star$, for each
$E,F \in \boldsymbol{\overline{F}}(D)$.  A spectral semistar operation
is stable \cite[Lemma 4.1]{FH2000}.

If $\star$ is a semistar operation on $D$, we denote by
${\tilde{\star}}$ the semistar operation
$\star_{\mathcal{M}(\star_{\!_f})}$ induced by the set
$\mathcal{M}(\star_{\!_f})$ of the quasi--$\star_{\!_f}$--maximal
ideals of $D$.  The semistar operation ${\tilde{\star}}$ is stable and
of finite type and ${\tilde{\star}} \leq \star_{\!_f}$ (cf.  \cite[p.
181]{FH2000}, where the semistar operation ${\tilde{\star}}$ is defined,
in an equivalent way, by using localizing systems, and also
\cite[Section 2]{AC2000} for an analogous construction in the star
setting).  Note that when $\star = v$ (where, as usual, $v$ denotes the
(semi)star operation defined by $E^v := (D:(D:E))$, for each $E\in
\boldsymbol{\overline{F}}(D)$), then ${\tilde{\star}}$ coincides with
the (semi)star operation denoted by $w$ by Wang Fanggui and R.L.
McCasland (cf.  \cite{WM97}, \cite{WM97b} and \cite{WM99}).

   The following lemma is not difficult to prove  (cf.  \cite[Corollary
 3.5(2)]{FL03} and, for the analogous result in case of star
 operations, \cite[Theorem 2.16]{AC2000}).

\begin{lemma} \label{lemma:max}
Let $\star$ be a semistar operation on an integral domain $D$.
Then, $\mathcal{M}(\star_{\!_f})=\mathcal{M}({\tilde{\star}})$. \hfill
$\Box$
\end{lemma}

In the  next proposition, we recall how a semistar operation on an
integral domain $D$ induces canonically a semistar operation on an
overring $T$ of $D$ (cf.  \cite[Lemma 45]{OM94}, and, for the notation
used here, \cite{FP03}).

\begin{prop} \label{stari}
 Let $D$ be an integral domain and $T$ an overring of
$D$. Let $\iota:D \hookrightarrow T$ be the embedding of $D$ in
$T$, and let $\star_{\iota} : \boldsymbol{\overline{F}}(T) \to
\boldsymbol{\overline{F}}(T)$ be defined by $E^{\star_{\iota}} :=
E^\star$.  Then:
\begin{enumerate}
\item[(1)] $\star_{\iota}$ is a semistar operation on $T$.
\item[(2)] If $\star$ is of finite type on $D$, then $\star_{\iota}$ is
of finite type on $T$.
\item[(3)] If $T=D^\star$, then $\star_{\iota}$ is a (semi)star
operation on  $D^\star$.
\item[(4)] If $\star$ is stable, then
$\star_{\iota}$ is stable.  \hfill $\Box$
\end{enumerate}
\end{prop}

If $R$ is a ring and $X$ an indeterminate over $R$, then the ring
$R(X):=\{f/g \, \vert \; f,g \in R[X] \mbox{ and }
\boldsymbol{c}(g)=R \}$ (where $\boldsymbol{c}(g)$ is the content
of the polynomial $g$) is called the \emph{ Nagata ring} of $R$
\cite[Proposition 33.1]{Gil92}.

The following result is proved in \cite[Proposition 3.1]{FL03}
(cf.  also \cite[Proposition 2.1]{K89}).

\begin{prop}  \label{prop:nagata}
Let $\star$ be a nontrivial semistar operation on an integral
domain $D$ and set $N(\star)  := N_D(\star):=\{h \in D[X] \, \vert \; h
\neq 0 \mbox{ and } (\boldsymbol{c}(h))^\star=D^\star \}$.  Then:
\begin{enumerate}
\item[(1)] $N(\star)$ is a saturated multiplicative subset of
$D[X]$ and $N(\star)=N(\star_{\!_f})=D[X] \smallsetminus \bigcup \{Q[X]
\, \vert \; Q \in \mathcal{M}(\star_{\!_f}) \}$.
\item[(2)] $\Max(D[X]_{N(\star)})= \{Q[X]_{N(\star)} \, \vert  \;   Q \in
\mathcal{M}(\star_{\!_f}) \}$ and $\mathcal{M}(\star_{\!_f})$ coincides with
the canonical image in\ $\Spec(D)$ of\
$\Max \left((D[X])_{N(\star)} \right)$.
\item[(3)]$D[X]_{N(\star)} =
\bigcap \{D_Q(X) \, \vert \; Q \in \mathcal{M}(\star_{\!_f})\}$.  \hfill
$\Box$
\end{enumerate}

\end{prop}

We set $\Na(D,\star):= D[X]_{N(\star)}$ and we call it \emph{the
Nagata ring of $D$ with respect to the semistar operation
$\star$}. Obviously, $\Na(D,\star)=\Na(D,\star_{\!_f})$ and, when
$\star=d$ (the  identity (semi)star operation) on $D$, then $\Na(D,d)=D(X)$.


\section{Semistar Invertibility}


Let $\star$ be a semistar operation on an integral domain $D$. Let
$I \in \boldsymbol{F}(D)$,  we say that $I$ is \emph{$\star$--invertible}
if $\left(II^{-1}\right)^\star=D^\star$.    In particular when $\star = d$
[respectively,\ $v$\,,\; $t\ (:=v_{_{\!f}})$\,,\; $w\ (:= \widetilde{v}$)\ ] is
the identity (semi)star operation [respectively,\ the $v$--operation,\; the
$t$--operation,\; the $w$--operation\ ] we reobtain the classical notion
of \it invertibility \rm [respectively,\ \it $v$--invertibility,\;
$t$--invertibility,\; $w$--invertibility\ \rm ] of a fractional ideal.

\begin{lemma} \label{lemma:inv}
Let $\star, {\star_1}, {\star_2}$ be semistar operations on an
integral domain $D$. Let  $\Inv(D, \star)$ be the set of all
$\star$--invertible fractional ideals of $D$ and $\Inv(D)$ (instead of
$\Inv(D, d)$) the set of all invertible fractional ideals of $D$.  Then:
\begin{enumerate}
\item[(0)] $D \in  \Inv(D, \star)$.

\item[(1)] If ${\star_1} \leq {\star_2}$, then $ \Inv(D, \star_1) \subseteq
\Inv(D, \star_2)$. In particular, $\Inv(D) \subseteq \Inv(D,
\tilde{\star}) \subseteq \Inv(D,\star_{\!_f}) \subseteq \Inv(D,
\star)$.

\item[(2)] $I,J \in  \Inv(D, \star)$ if and only if $IJ \in
 \Inv(D, \star)$.

\item[(3)] If $I \in  \Inv(D, \star)$ then $I^{-1} \in
 \Inv(D, \star)$..

\item[(4)] If $I \in  \Inv(D, \star)$ then $I^v \in
 \Inv(D, \star)$.

\end{enumerate}
\end{lemma}
\begin{proof}
(0) and (1) are obvious.

(2) Note that, if $I,J \in   \Inv(D, \star)  $, then $D^\star =
\left(II^{-1}\right)^\star\left(JJ^{-1}\right)^\star \subseteq
\left(II^{-1}JJ^{-1}\right)^\star \subseteq \left(IJ(IJ)^{-1}\right)^\star
\subseteq D^\star$.  Thus, $IJ \in  \Inv(D, \star)  $.  Conversely, if $IJ \in
  \Inv(D, \star)$, then
 $D^\star=((IJ)(D:IJ))^\star=(I(J(D:IJ)))^\star$.  Since $(J(D:IJ))
 \subseteq (D:I)$, it follows $(I(D:I))^\star=D^\star$.  Similarly,
 $(J(D:J))^\star=D^\star$.

(3) $D^\star = \left(II^{-1}\right)^\star \subseteq
\left((I^{-1})^{-1}I^{-1}\right)^\star \subseteq D^\star$.

(4) follows from (3).
\end{proof}

\vskip 4pt

\begin{rem} \label{rem:2.2}
      \bf (a) \rm Note that $D$ is the unit element of $  \Inv(D, \star)  $
      with respect to the \emph{standard multiplication} of fractional ideals
      of $D$.  Nevertheless, $  \Inv(D, \star)  $ is \emph{not} a group in
      general (under the standard multiplication), because for $I \in
       \Inv(D, \star)  $, then $I^{-1} \in  \Inv(D, \star)  $, but $II^{-1} \neq D$, if $I
      \not \in \Inv (D)$.  For instance, let $k$ be a field, $X$ and $Y$ two
      indeterminates over $k$, and let $D:=k[X,Y]_{(X,Y)}$.  Then $D$ is a
      local Krull domain, with maximal ideal $M:=(X,Y)D$.  Let $\star=v$, then
      clearly $M^v=D$, since $\mbox{ht}(M)=2$, thus $M$ is $v$--invertible but
      $M$ is not invertible in $D$, since it is not principal.  Therefore
      $(MM^{-1})^v=D$, but $M= MM^{-1} \subsetneq D$.  We will discuss
      later what happens if we consider the semistar (fractional) ideals
      semistar invertible with the ``semistar product''.

 \bf (b) \rm Let $I \in \boldsymbol{F}(D)$.  Assume that $I \in
 \Inv(D, \star)$ and $(D^\star: I) \in \boldsymbol{F}(D)$, then  we will
 see later that $(D^\star :I) = (D:I)^\star$ (Lemma
 \ref{lemma:quasistarinv}, Remark \ref{rem:2.12}(d1) and Proposition
 \ref{pr:2.16}), more precisely that:

 \centerline{$\left(I^{-1}\right)^\star = \left(D:I^v \right)^\star =
\left(D^\star :I\right)^\star = (D^\star :I) = \left(I^\star\right)^{-1}\,.$}

 However, in this situation, we may not conclude that $(D^\star :I)$ (or,
 $(D:I)^\star$) belongs to $\Inv(D, \star)$ (even if $(D:I)\in \Inv(D,
 \star)$, by Lemma \ref{lemma:inv}(3)).  As a matter of fact, more generally, if $J \in
 \Inv(D, \star)$ and $J^\star \in \boldsymbol{F}(D)$, then $J^\star$ does
 not belong necessarily to $\Inv(D, \star)$.

 For instance, let $K$ be a
 field and $X, Y$ two indeterminates over $K$, set $T: =K[X, Y]$ and $D := K
 + YK[X, Y]$.  Let $\star_{\{T\}}$ be the semistar operation on $D$
 defined by $E^{\star_{\{T\}}} := ET$, for each $E\in
 \boldsymbol{\overline{F}}(D)$.  Then $J:= YD$ is obviously invertible
 (hence  $\star_{\{T\}}$--invertible) in $D$ and $J^{\star_{\{T\}}}= JT =
 YT= YK[X, Y] = (D:T)$ is a nonzero (maximal) ideal of $D$ (and, at the same
 time, a (prime) ideal of $T$),\ but $J^{\star_{\{T\}}}$ is not
 $\star_{\{T\}}$--invertible in $D$, because
 $\left(J^{\star_{\{T\}}}(D:J^{\star_{\{T\}}})\right)^{\star_{\{T\}}}
 =(JT(D:JT))T = (YT(D:YT))T= (YTY^{-1}(D:T))T = (T(YT))T =YT \subsetneq T
 = D^{\star_{\{T\}}}$.
%
%

\bf (c) \rm Note that the converses of (3) and (4) of Lemma
\ref{lemma:inv} are  not true in general.  For instance, take an
integral domain $D$ that is not an H--domain (recall that an \it
H--domain \rm  is an integral domain $D$ such that, if $I$ is an
ideal of $D$ with $I^{-1}=D$, then there exists a finitely generated $J
\subseteq I$, such that $J^{-1}=D$ \cite[Section 3]{GV77}).  Then, there
exists an ideal $I$ of $D$ such that $I^v=I^{-1}=D$ and $I^t \subsetneq
D$.  It follows that $\left(I^{-1}I^v \right)^t=D$ ( and so, $I^{-1}$
and $I^v$ are $t$--invertibles), but $\left(II^{-1}\right)^t=I^t
\subsetneq D$, that is, $I$ is not $t$--invertible.  On the other hand,
note that, trivially, $I$ is $v$--invertible.

An explicit example is given by a $1$--dimensional non discrete
valuation domain $V$ with maximal ideal $M$. Clearly, $V$ is not
an H--domain \cite[(3.2d)]{GV77}, $M^{-1}=M^v=V$ \cite[Exercise
12 p.431]{Gil92} and $M^t = \bigcup \{J^v \vert J \subseteq M, J
\mbox{ finitely generated} \} = \bigcup \{J \vert J \subseteq M, J
\mbox{ finitely generated} \} = M \subsetneq V$. In this case,
$M^{-1}$ and $M^{v}$ are obviously $t$--invertibles, but $M$ is not
$t$--invertible.
\end{rem}

\vskip 8pt If $I \in \boldsymbol{\overline{F}}(D)$, we say that
$I$ is \emph{$\star$--finite} if there exists $J \in
\boldsymbol{f}(D)$ such that $J^\star=I^\star$.  It is immediate
to see that if ${\star_1} \leq {\star_2}$ are semistar operations
and $I$ is ${\star_1}$--finite, then $I$ is ${\star_2}$--finite.
In particular, if $I$ is $\star_{\!_f}$--finite, then it is
$\star$--finite.

 We notice that, in the previous definition of $\star$--finite, we
do not require that $J \subseteq I$. Next result shows that, with
this ``extra'' assumption, $\star$--finite is equivalent to
$\star_{\!_f}$--finite.

\begin{lemma} \label{lm:2.4}
    Let $\star$ be a semistar operation on an integral
domain $D$ with quotient field $K$. Let $I \in
\boldsymbol{\overline{F}}(D)$. Then, the following are equivalent:
\begin{enumerate}
\item[(i)] $I$ is $\star_{\!_f}$--finite.
\item[(ii)] There exists $J \subseteq I$, $J \in \boldsymbol{f}(D)$ such that
$J^\star = I^\star$.
\end{enumerate}
\end{lemma}
\vskip -4pt \begin{proof} It is clear that (ii) implies (i), since
$J^\star=J^{\star_{\!_f}}$, if $J$ is finitely generated. On the other
hand, suppose $I$ $\star_{\!_f}$--finite. Then,
$I^{\star_{\!_f}}=J_0^{\star_{\!_f}}$, with $J_0=(a_1,a_2, \ldots, a_n)D$,
for some family $\{a_1,a_2, \ldots, a_n\} \subseteq K$. Since $J_0
\subseteq I^{\star_{\!_f}}$, there  exists a finite family of finitely
generated fractional ideals of $D$, $J_1, J_2, \ldots, J_n \subseteq
I$, such that $a_i \in J_i ^{\star}$, for $i=1,2, \ldots, n$.  It
follows that $I^{\star_{\!_f}} = J_0 ^{\star_{\!_f}} \subseteq
\left(J_1^{\star_{\!_f}}+ J_2^{\star_{\!_f}} + \ldots +
J_n^{\star_{\!_f}}\right)^{\star_{\!_f}} = (J_1 + J_2 + \ldots +
J_n)^{\star_{\!_f}} \subseteq I^{\star_{\!_f}}$.  Set $J:=J_1 + J_2 +
\ldots + J_n$.  Then, $J$ is finitely generated, $J \subseteq I$ and
$J^{\star_{\!_f}} = I^{\star_{\!_f}}$, thus $J^\star = I^\star$.
\end{proof}

\begin{rem} Extending the terminology introduced by Zafrullah in the
star setting \cite{Z1} (cf.  also \cite[p.  433]{Z2}),
given a semistar operation on an integral domain $D$, we can say that
 $I\in \boldsymbol{\overline{F}}(D)$ is  \it strictly $\star$--finite
 \rm if $I^\star = J^\star$, for some $J \in \boldsymbol{{f}}(D)$, with
 $J \subseteq I$.  With this terminology, Lemma \ref{lm:2.4} shows that
 $\star_{\!_f}$--finite coincides with strictly $\star_{\!_f}$--finite.  This
 result was already proved, in the star setting, by Zafrullah
 \cite[Theorem 1.1]{Z1}. \ Note  that Querr\'e  studied the strictly
 $v$--finite ideals \cite{Q}, using often the terminology of  \it
 quasi--finite ideals.\rm

For examples of  $\star$--finite ideals that are not
$\star_{\!_f}$--finite (when $\star$ is the $v$--operation),  see
\cite[Section (4c)]{GH97}, where are described domains with all
the ideals $v$--finite (called DVF--domains), that are not Mori domains
(that is, such that not all the ideals are $t$--finite).  \end{rem}



\begin{lemma} \label{lm:2.6}

    Let $\star$ be a semistar operation on an integral
domain $D$ and let $I \in \boldsymbol{F}(D)$.  Then $I$ is
$\star_{\!_f}$--invertible if, and only if, $\left(I^\prime I^{\prime
\prime}\right)^\star=D^\star$, for some $I^{\prime} \subseteq I, I^{\prime
\prime} \subseteq I^{-1}$, and $I^\prime, I^{\prime \prime} \in
\boldsymbol{f}(D)$.  Moreover, ${I^\prime}^{\star}=I^\star$ and
${I^{\prime \prime}}^\star=\left(I^{-1}\right)^\star$.
\end{lemma}
\begin{proof}
The `` if'' part is trivial.  For the ``only if'': if
$\left(II^{-1}\right)^{\star_{\!_f}}=D^{\star_{\!_f}}$, then $H^\star =
D^\star$ for some $H \subseteq II^{-1}$, $H \in \boldsymbol{f}(D)$.
Therefore, $H=(h_1, h_2 \ldots , h_n)D$, with $h_i = x_{1,i}y_{1,i} +
x_{2,i}y_{2,i} + \ldots + x_{k_i,i}y_{k_i,i}$, with the $x$'s in $I$ and
the $y$'s in
$I^{-1}$. 
Let $I^{\prime}$ be the  (fractional)  ideal of $D$ generated by the $x$'s
and let $I^{\prime \prime}$ be the (fractional) ideal of $D$ generated
by the $y$'s.  Then, $H \subseteq I^{\prime}I^{\prime \prime} \subseteq
II^{-1}$ and so $D^\star = \left(I^{\prime}I^{\prime
\prime}\right)^\star$, and, thus, also $D^\star=\left(I^\prime
I^{-1}\right)^\star=\left(II^{\prime \prime}\right)^\star$.  Moreover,
$I^\star = \left(ID^\star\right)^\star= \left(I\left(I^\prime
I^{-1}\right)^\star\right)^\star= \left((II^{-1})^\star
I^{\prime}\right)^\star=\left(D^\star
I^\prime\right)^\star={I^\prime}^\star$.  In a similar way, we obtain
also that ${I^{\prime \prime}}^\star= \left(I^{-1}\right)^\star$.
\end{proof}

  A classical result due to Krull \cite[Th\'eor\`eme 8, Ch.  I, \S
4]{Jaffard60} shows that, for a star operation of finite type,
star--invertibility implies star--finiteness.  The following result gives
a more complete picture of the situation in the general semistar
setting.

\begin{prop} \label{cor:invfin}
Let $\star$ be a semistar operation on an integral domain $D$. Let
$I \in \boldsymbol{F}(D)$.    Then $I$ is $\star_{\!_f}$--invertible
if and only if $I$ and $I^{-1}$ are $\star_{\!_f}$--finite (hence, in
particular, $\star$--finite) and $I$ is $\star$--invertible.
\end{prop}
  \begin{proof}
The  ``only if'' part follows from Lemma \ref{lm:2.6} and from the
fact that $\star_{\!_f}\leq \star$.

For the ``if'' part, note that by assumption $I^{\star_{\!_f}}
={J'}^{\star_{\!_f}}= {J'}^\star$ and $(I^{-1})^{\star_{\!_f}}
={J''}^{\star_{\!_f}}= {J''} ^\star$, with $J', J'' \in
\boldsymbol{f}(D)$.  Therefore:
$$
 \left(II^{-1}\right)^{\star_{\!_f}}= \left(J'J''\right)^{\star_{\!_f}}
 =\left(J'J''\right)^\star = \left({J'}^\star {J''}^\star\right)^\star =
 (I^\star\left(I^{-1})^\star\right)^\star = \left(II^{-1}\right)^\star = D^{\star}\,.
 $$
 \end{proof}

  Next goal is to investigate when the $\star$--inver\-ti\-bi\-li\-ty coincides
with the $\star_{\!_f}$--inver\-ti\-bi\-li\-ty.

Let $\star$ be a semistar operation on an integral domain $D$, we
say that $D$ is an \it $H(\star)$--domain \rm if, for each nonzero
integral ideal $I$ of $D$ such that $I^\star = D^\star$, there
exists $J \in \boldsymbol{{f}}(D)$ with $J \subseteq I$ and
$J^\star =D^\star$.  It is easy to see that,  for $\star =v$, the
H$(v)$--domains coincide with the H--domains introduced by Glaz
and Vasconcelos (Remark  \ref{rem:2.2}(c)).

\begin{lemma} \label{le:2.8} Let $\star$ be a semistar operation
on an integral domain $D$.  Then $D$ is an H$(\star)$--domain if and
only if each quasi--$\star_{_{\!f}}$--maximal ideal of $D$ is a
quasi--$\star$--ideal of $D$.
        \end{lemma}
     \begin{proof} Assume that $D$ is an H$(\star)$--domain.  Let $Q =
     Q^{\star_{_{\!f}}} \cap D$ be a quasi--$\star_{_{\!f}}$--maximal ideal
     of $D$.  Assume that $Q^\star = D^\star$.  Then, for some $J \in
     \boldsymbol{f}(D)$, with $J \subseteq Q$, we have $J^\star = D^\star$,
     thus $Q^{\star_{_{\!f}}}=D^\star$, which leads to a contradiction.
     Therefore $Q^{\star_{_{\!f}}} \cap D \subseteq Q^\star \cap D \subsetneq
     D$ and, hence, there exists a quasi--$\star_{_{\!f}}$--maximal ideal of
     $D$ containing $Q^\star \cap D$.  This is possible only if
     $Q^{\star_{_{\!f}}} \cap D = Q^\star \cap D$.

      Conversely, let $I$ be a nonzero ideal of $D$ with the property $I^\star
      = D^\star$.  Then, necessarily $I \not\subseteq Q$ for each
      quasi--$\star_{_{\!f}}$--maximal ideal of $D$ (because, otherwise, by
      assumption $I\subseteq Q = Q^{\star_{_{\!f}}} \cap D = Q^\star \cap D$,
      and so $I^\star \subseteq Q^\star \subsetneq D^\star$).  Therefore
      $I^{\star_{_{\!f}}}= D^\star$.
      \end{proof}

        Next result provides several characterizations of the
      H$(\star)$--domains;\ note that, in the particular case that $\star = v$,
      the equivalence (i) $\Leftrightarrow$ (iii) was already known
      \cite[Proposition 2.4]{HZ88} and the equivalence (i) $\Leftrightarrow$
      (iv) was considered in \cite[Proposition 5.7]{WM97}.

\begin{prop} \label{prop:2.9}  Let $\star$ be a semistar operation
on an integral domain $D$. The following are equivalent:
\begin{enumerate}
\item[{\rm (i)}] $D$ is an $H(\star)$--domain;
\item[{\rm (ii)}] for each $I \in\boldsymbol{{F}}(D)$,\ $I$ is
$\star$--invertible if and  only if $I$ is $\star_{\!_f}$--invertible;
\item[{\rm (iii)}] $\mathcal{M}(\star_{\!_f}) = \mathcal{M}(\star)$;
\item[{\rm (iv)}] $\mathcal{M}(\widetilde{\star}) = \mathcal{M}(\star)$.
\end{enumerate}
\end{prop}
\begin{proof} Obviously, (iii)\ $\Leftrightarrow$\ (iv) by Lemma \ref{lemma:max} and
(iii)\ $\Leftrightarrow$\ (i) by Lemma \ref{le:2.8}, recalling that a
quasi--$\star$--ideal is also a quasi--$\star_{\!_f}$--ideal.

  (iii)\ $\Rightarrow$\ (ii). Let $I$ be a $\star$--invertible ideal of
$D$. Assume that $I$ is not $\star_{\!_f}$-invertible. Then, there
exists a quasi--$\star_{\!_f}$--maximal ideal $M$ such that $II^{-1}
\subseteq M$. But $M$ is also quasi-$\star$-maximal, since
$\mathcal{M}(\star_{\!_f})=\mathcal{M}(\star) $. Thus $M^\star
\subsetneq D^\star$. It follows that $(II^{-1})^\star \subseteq
M^\star \subsetneq D^\star$, a contradiction.  Hence $I$ is
$\star_{\!_f}$-invertible.


(ii)\ $\Rightarrow$\ (i)  Let $I$ be a  nonzero integral ideal $I$ of $D$ such that $I^\star =
D^\star$. Then, $I \subseteq II^{{-1}} \subseteq D$ implies that
$\left(II^{{-1}}\right)^{\star} = D^{\star}$, that is $I$ is
$\star$--invertible.  By assumption, it follows that  $I$ is
$\star_{\!_f}$--invertible, and so $I$ is
$\star_{\!_f}$--finite (Proposition \ref{cor:invfin}). By Lemma
\ref{lm:2.4}, we conclude that there exists $J \in \boldsymbol{{f}}(D)$ with $J \subseteq I$
and $J^\star  =I^{\star}=D^{\star}$.
 \end{proof}
 \smallskip

   Let $\star$ be a semistar operation of $D$. If we denote by
$\iota: D \hookrightarrow D^\star$ the embedding of $D$ in
$D^\star$ and by $\star_{\iota}$ the (semi)star operation
canonically induced on $D^\star$ by $\star$ (defined as in
Proposition \ref{stari}), we note that, if $I \in
 \Inv(D, \star)$, then $I^\star \in  \Inv(D^\star,
\star_{\iota})$. As a matter of fact, we have: $D^\star =
\left(II^{-1}\right)^\star=\left(I^\star(D:I)^\star\right)^\star
\subseteq \left(I^\star(D^\star:I^\star)\right)^\star =
\left(I^\star(D^\star:I^\star)\right)^{\star_{\iota}}\subseteq
\left(D^\star\right)^\star = D^\star$.

Next example shows that the converse does not hold (in other words
$I^\star$ may be in $\Inv(D^\star, \star_{\iota})$, with $I\in
\boldsymbol{F}(D) \smallsetminus  \Inv(D, \star)  $), even if
$\star$ is a semistar operation stable and of finite type.

\begin{ex} \label{ex:almded}
Let $D$ be an almost Dedekind domain, that is not a Dedekind
domain (cf.  for instance \cite[Section 2 and the references]{Gil90}).
Then, in $D$ there exists a prime (= maximal) ideal $P$, such that
$P$ is not invertible (otherwise, $D$ would be a Dedekind domain).
Then, $P^{-1}=D$ \cite[Corollary 3.1.3]{FHP97}, since $D$ is a Pr\"ufer
domain.  Consider the semistar operation $\star := \star_{\{P\}}$.  Let
$\iota_P: D \hookrightarrow D_P$ be the canonical embedding, then
$P^\star= PD_P$ is $\star_{\iota_P}$--invertible, since $D_P$ is a DVR,
but $\left(PP^{-1}\right)^\star=(PD)^\star=P^\star= PD_P \subsetneq
D_P=D^\star$, thus $P$ is not $\star$--invertible.
\end{ex}

Let $\iota:D \hookrightarrow D^\star$ be the canonical embedding,
then, we say that $I \in \boldsymbol{\overline{F}}(D)$ is
\emph{quasi--$\star$--invertible} if $I^\star \in
 \Inv(D^\star, \star_{\iota})$ (that is, if
$\left(I(D^\star:I)\right)^\star=D^\star)$.  Note that $I^\star \in
 \Inv(D^\star, \star_{\iota})$ implies that $I^\star \in
\boldsymbol{F}(D^\star)$.  We denote by $  \QInv(D, \star)  $ the set of all
quasi--$\star$--invertible $D$--submodules of $K$ and, when $\star
=d$, we set $\QInv(D)$, instead of $\QInv(D, d)$. We have already
noticed that $  \Inv(D, \star)  \subseteq  \QInv(D, \star)
$ and that the inclusion can be proper.   Moreover, it is obvious that
$\QInv(D) = \Inv(D)$.

We have the following straightforward necessary and sufficient
condition for a $D$--submodule of $K$ to be
quasi--$\star$--invertible.

\begin{lemma} \label{lemma:quasistarinv}
Let $\star$ be a semistar operation on an integral domain $D$ and
$I \in   \boldsymbol{\overline{F}}(D)$.  Then, $I$ is
quasi--$\star$--invertible if and only if there exists $H \in
\boldsymbol{\overline{F}}(D)$ such that $(IH)^\star=D^\star$.  \hfill
$\Box$
\end{lemma}

Next we prove an analogue of Lemma \ref{lemma:inv} for quasi--
$\star$--invertible ideals.

\begin{lemma} \label{lemma:qinv}
Let $\star, {\star_1}, {\star_2}$ be semistar operations on an
integral domain $D$. Then:
\begin{enumerate}
\item[(0)] $D^\star \in   \QInv(D, \star)  $.
\item[(1)] If ${\star_1} \leq {\star_2}$, then $\QInv(D, {\star_1})
\subseteq \QInv (D, {\star_2})$.\ In particular, we have

\centerline{$\QInv(D) \subseteq  \QInv(D, \tilde{\star}) \subseteq
\QInv (D, \star_{\!_f}) \subseteq  \QInv(D, \star) $.  }

\item[(2)]
$I,J \in \QInv(D,\star)$ if and only if $IJ \in \QInv(D,
\star)$. \item[(3)] If $I \in  \QInv(D, \star)$,  then
$(D^\star:I) \in  \QInv(D, \star)$.   \item[(4)] If $I \in
\QInv(D, \star)  $, then $I^{v(D^\star)}:= (D^\star:(D^\star:I)) \in
 \QInv(D, \star)$.
\end{enumerate}
\end{lemma}
\begin{proof} (0) and (1) are straightforward. To prove  (2) we notice that
$I, J \in \!  \QInv(D, \star)$ if and only if $I^\star,J^\star \!\in\!
\Inv(D^\star, \star_{\iota})$, where $\star_{\iota}$ is defined as
above.  It follows (from Lemma \ref{lemma:inv}(2)) that $I,J \in
\QInv(D, \star)  $ if and only if $I^\star J^\star \in
\Inv(D^\star, \star_{\iota})$. It is easy to see that this happens if
and only if $(IJ)^\star \in \Inv(D^\star, \star_{\iota})$, that is, if
and only if $IJ \in  \QInv(D, \star)  $.  (3)  is clear and
(4)   is an immediate consequence of Lemma \ref{lemma:inv}(4)  and
of the fact that $\left(v(D^\star)\right)_{\iota} = v_{D^\star}$,
where $v_{D^\star}$ is the $v$--operation of $D^\star$, $\iota$ is the
canonical embedding of $D$ in $D^\star$ and $v(D^\star)$ is the
semistar operation on $D$, defined by $E^{v(D^\star)} :=
(D^{\star}:(D^{\star} :E))$, for each $E \in
\boldsymbol{\overline{F}}(D)$ \ (note that, obviously, $\star \leq
v(D^\star)$).
\end{proof}

\begin{cor} \label{cor:2.10}
Let $\star$ be a semistar operation on an integral domain $D$, let
$v(D^\star)$ be the semistar operation on $D$, defined in the proof of
Lemma \ref{lemma:qinv}(4) and let $I \in \boldsymbol{\overline{F}}(D)$.
If $I$ is quasi--$\star$--invertible, then $I$ is
quasi--$v(D^\star)$--invertible and $I^\star=I^{v(D^\star)}$.
\end{cor}
  \begin{proof}  Let $\iota$ be the canonical
embedding of $D$ in $D^\star$. As we noted in the proof of Lemma
\ref{lemma:qinv} (4),   $\left(v(D^\star)\right)_{\iota} =
v_{D^\star}$.  Then, in order to show that $I^\star$ is
quasi--$v(D^\star)$--invertible, we prove that $I^\star$ is
$v_{D^\star}$--invertible.  But $\star_\iota$ is a (semi)star operation
on $D^\star$ and $I^\star$ is $\star_\iota$--invertible, then (Lemma
\ref{lemma:inv} (1)) $I^\star$ is $v_{D^\star}$--invertible, since
$\star_{\iota} \leq v_{D^\star}$ \cite[Theorem 34.1(4)]{Gil92}.
Therefore $I$ is quasi--$v(D^\star)$--invertible and $I^{\star}=
\left(I^{v(D^\star)}\right)^{\star}=I^{v(D^\star)}$, since $
(D^{\star}:I) = \left(D^{\star}:I^{v(D^\star)}\right)$ (cf.  also
\cite[p.  433]{Z2} or \cite[Lemma 2.1(3)]{CP03}, and Remark
\ref{rem:2.12}(b1)).  \end{proof}

\begin{rem}\label{rem:2.12}
    \bf (a) \rm Note that if $I$ is a quasi--$\star$--invertible \sl ideal
    \rm of $D$, then every ideal $J$ of $D$, with $I \subseteq J \subseteq
    I^\star \cap D$, is also quasi--$\star$--invertible.

      More precisely, \it  let
    $I, J \in\boldsymbol{{F}}(D)$ [respectively, $I, J
    \in\boldsymbol{\overline{F}}(D)$], assume that $J\subseteq I$,\ $J^\star
    =I^\star$ and that $I$ is $\star$--invertible [respectively,
    quasi--$\star$--invertible] then $J$ is $\star$--invertible [respectively,
    quasi--$\star$--invertible].  \rm

    Conversely, \it let $I, J
    \in\boldsymbol{\overline{F}}(D)$, assume that $J\subseteq I$,\ $J^\star =I^\star$ and that
    $J$ is quasi--$\star$--invertible then $I$ is quasi--$\star$--invertible
    (but not necessarily $\star$--invertible, even if $J$ is
    $\star$--invertible).\rm

    As a matter of fact, if $I$ is $\star$--invertible, then $D^\star =
    (I(D:I))^\star = (J(D:I))^\star \subseteq (J(D:J))^\star\subseteq D^\star$.
    The quasi--$\star$--invertible case is similar.  Conversely, if $J$ is
    quasi--$\star$--invertible then
    $D^\star =
    \left(J(D^\star:J)\right)^\star = \left(I(D^\star :J)\right)^\star $, thus
    $I$ is quasi--$\star$--invertible and  $(D^\star :J) =\left(D^\star
    :J\right)^\star =\left(D^\star :I \right)^\star=(D^\star :I)$ (cf.  also
    (d1)).

    Example \ref{ex:almded} shows the parenthetical part of the statement.
    Let $D$, $P$ and $\star$ be as in Example \ref{ex:almded}.
     Note that $P^{\star}$ is principal in (the
    DVR) $D^{\star}=D_{P}$, thus $P^{\star}= PD_P =tD_{P}$, for some nonzero
    $t\in PD_P$.  Therefore, if $J:= tD$, then $J^\star = P^\star$, i.e. $P$
    is $\star$--finite.  We already observed that $P$ is
    quasi--$\star$--invertible but not $\star$--invertible, even if
    obviously $J$ is ($\star$--)invertible.

    \bf (b) \rm Let $I,\ H',\ H'',\ J,\ L \in \boldsymbol{\overline{F}}(D)$.
    The following properties are straightforward:
    \begin{enumerate}
        \item[(b1)] $\left(IH'\right)^\star=D^\star= \left(IH''\right)^\star \;
        \; \Rightarrow\; \; {H'}^\star = {H''}^\star = \left(D^\star :
        I\right)^\star =(D^\star : I)\,.$
        \item[(b2)] $I \in   \QInv(D, \star)  $, \,
        $IJ\subseteq IL \; \; \Rightarrow\; \; {J}^\star \subseteq {L}^\star\,.$
        \item[(b3)] $I \in   \QInv(D, \star)  $, \, $J\subseteq I^\star \; \;
        \Rightarrow\; \; \exists\ L \in \boldsymbol{\overline{F}}(D), \
        (IL)^\star = J^{\star}\,.$\newline [Take $L:= (D^\star :I)J$.\ ]
        \item[(b4)] $I,\ J \in   \QInv(D, \star)  $, \, $(IL)^\star = J^{\star} \; \;
        \Rightarrow\; \; L\in   \QInv(D, \star)  \,.$
        \newline
        [Set $H := I(D^\star:J)$, and note that $(LH)^\star = D^\star$.\ ]
        \item[(b5)] $I,\
        J \in   \QInv(D, \star)   \; \; \Rightarrow\; \; (D^\star :
        IJ)=\left(D^\star : IJ\right)^\star = \left(\left(D^\star
        : I\right) \left(D^\star : J\right)\right)^\star\,.$
        \item[(b6)] $I,\ J
        \in   \QInv(D, \star)   \; \; \Rightarrow\; \; \exists\ L\in   \QInv(D, \star)  ,\
        L \subseteq I^\star, \ L \subseteq J^\star\,.$ \newline [Take $z\in K$, $z
        \neq 0$, such that $zI \subseteq D^\star,\ zJ \subseteq D^\star,$ and
        set $L:= zIJ$.\ ]
        \item[(b7)] $I,\ J\in   \QInv(D, \star) \,, \, \, \, I+J
        \in   \QInv(D, \star)   \,\, \Rightarrow\,\,I^{v(D^\star)} \cap
        J^{v(D^\star)}\in   \QInv(D, \star)\,.$
         \newline
         [Recall that $\star \leq v(D^\star)$ and note that: \newline
         $\left(\left(D^\star : I \right)\!\left(D^\star :
         J\right)\!(I+J)\right)^\star \!\!  = \!  \left( \left( \left(D^\star :
         I\right)I\right)^\star\!  \left(D^\star : J\right)\!  + \!
         \left(D^\star : I\right) \!\left(\left(D^\star : J\right)J\right)^\star
         \right)^\star$ \!\!  = \!$\left(\left(D^\star : J\right)+ \left(D^\star
         : I\right)\right)^\star $ = $\left( \left(D^\star : J^{v(D^\star)}
         \right) + \left(D^\star : I^{v(D^\star)}\right) \right)^\star $ $\,
         \Rightarrow \,$ \newline $\left(\left(D^\star : I \right)\left(D^\star :
         J\right)(I+J)\right)^{v(D^\star) }$ = $\left( \left(D^\star :
         I^{v(D^\star)} \right) + \left(D^\star : J^{v(D^\star)}\right)
         \right)^{v(D^\star)}$ $\, \Rightarrow \,$ \newline
         $\left(D^{\star}:\left(\left(D^\star : I \right)\left(D^\star :
         J\right)(I+J)\right)\right)$ = $\left(D^{\star}:\left( \left(D^\star :
         I^{v(D^\star)} \right) + \left(D^\star : J^{v(D^\star)}\right)
         \right)\right) $ = $\left(D^{\star}: \left(D^\star : I^{v(D^\star)}
         \right)\right) \cap \left(D^{\star}: \left(D^\star : J^{v(D^\star)}
         \right)\right)$ = $ I^{v(D^\star)} \cap J^{v(D^\star)}$.\ ]

         \item[(b8)] $I, J\in
           \QInv(D, \star)\,, \,\, I^{v(D^\star)} \cap J^{v(D^\star)} \in   \QInv(D, \star)
        \,\,\Rightarrow\,\,
         I+J\in \QInv(D, v(D^\star)).$
         [Since $I^{v(D^\star)} \cap
         J^{v(D^\star)} = \left(D^{\star}:\left(\left(D^\star : I
         \right)\left(D^\star : J\right)(I+J)\right)\right)$ \ and hence \newline
         $\left(D^{\star}: \left(I^{v(D^\star)} \cap
         J^{v(D^\star)}\right)\right)= \left(\left(D^\star : I
         \right)\left(D^\star : J\right)(I+J)\right)^{v(D^\star)} $, then apply
         (b4) to conclude that $ I+J \in \QInv(D, v(D^\star))$.\ ]
        \end{enumerate}

         \bf (c) \rm A statement analogous to Corollary \ref{cor:2.10} holds for
         $\star$--invertibles:\ \it Let $\star$ be semistar operation on an
         integral domain $D$, let $v(D^\star)$ be the semistar operation on $D$,
         defined in the proof of Lemma \ref{lemma:qinv}(4) and let $I \in
         \boldsymbol{{F}}(D)$.  If $I$ is $\star$--invertible, then $I$ is
         $v(D^\star)$--invertible and $I^\star=I^{v(D^\star)}$.\rm

           \bf (d) \rm\sl Mutatis mutandis, \rm the statements in (b) hold for
$\star$--invertibles.  More precisely: Let $\star$ be a semistar
operation on an integral domain $D$ and let $I,\ H',\ H'',\ J,\ $ $ L \in
\boldsymbol{{F}}(D)$, then:
    \begin{enumerate}
        \item[(d1)] $I \in \Inv(D, \star)$,\, $(IH')^\star=D^\star=
        (IH'')^\star \; \; \Rightarrow \; \; {H'}^\star = {H''}^\star =
        \left(I^{-1}\right)^\star\,.
        $ \item[(d2)] $I \in  \Inv(D, \star)$,\,
        $IJ\subseteq IL \; \; \Rightarrow \; \; {J}^\star \subseteq
        {L}^\star\,.$
        \item[(d3)] $I \in  \Inv(D, \star)$,\,  $J\subseteq I^\star
        \; \; \Rightarrow \; \; \exists\ L \in \boldsymbol{{F}}(D), \ (IL)^\star
        = J^{\star}\,.$
        \item[(d4)] $I,\ J \in  \Inv(D, \star)$,\, $(IL)^\star =
        J^{\star} \; \; \Rightarrow \; \; L\in  \QInv(D, \star)\,\!, \;
        (D^\star :L)= \left(I(D:J)\right)^\star$\,.  \newline Note that, under
        the present hypotheses, $L\in  \Inv(D, \star)$ if and only if $ (D
        :L)^\star= \left(I(D:J)\right)^\star$\,.
        \item[(d5)]
        $I,\ J \in   \Inv(D, \star)  \; \; \Rightarrow \; \; \left(D:
        IJ\right)^\star = \left(\left(D : I\right) \left(D : J\right)\right)^\star\,.$
        \item[(d6)] $I,\ J \in   \Inv(D, \star)  \; \;
        \Rightarrow \; \; \exists\ L\in \Inv(D,\star),\  L \subseteq I, \ L
        \subseteq J\,.$
        \item[(d7)] $I,\ J\in \Inv(D,\star)\,, \, \; \; I+J \in  \Inv(D,
        \star)  \; \; \Rightarrow \; \; I^{v(D^\star)} \cap J^{v(D^\star)}\in
         \Inv(D, \star)  \,.$
        \item[(d8)] $I,\ J\in \Inv(D, \star)\,, \, \; \; I^{v(D^\star)}
        \cap J^{v(D^\star)} \in  \Inv(D, \star)  \; \; \Rightarrow \; \;
        I+J\in \Inv(D, v(D^\star))\,.$
      \end{enumerate}
    \end{rem}


  Our next goal is to extend Proposition \ref{cor:invfin} to the case
of quasi--$\star_{\!_f}$--invertibles.
We need the following:

\begin{lemma} \label{lemma:starifinite}
Let $\star$ be a semistar operation on an integral domain $D$ with
quotient field $K$, let  $\iota:D \hookrightarrow D^\star$ the
embedding of $D$ in $D^\star$, let $\star_{\iota}$ denote the
(semi)star operation canonically induced on $D^\star$ by $\star$
and let $I \in \boldsymbol{\overline{F}}(D)$. Then, $I$ is
$\star$--finite if and only if $I^\star$ is $\star_\iota$--finite.
\end{lemma}
\begin{proof}
If $I$ is $\star$--finite, then there exists $J \in
\boldsymbol{f}(D)$ such that $I^\star = J^\star$. It is clear that
$\left(JD^\star\right)^{\star_{\iota}} = I^{\star}$, with $JD^\star \in
\boldsymbol{f}(D^\star)$.  Thus, $I^\star$ is $\star_{\iota}$--finite.
Conversely, let $I^\star$ be $\star_{\iota}$--finite.  Then, there
exists $J_0 \in \boldsymbol{f}(D^\star)$, $J_0=(a_1, a_2, \ldots,
a_n)D^\star$, with $\{a_1,a_2, \ldots , a_n \} \subseteq K$, such
that $J_0 ^\star = J_0^{\star_\iota}=I^{\star_\iota}= I^\star$.  Set
$J=(a_1, a_2, \ldots, a_n)D \in \boldsymbol{f}(D)$. Then, $J^\star =
(a_1 D + a_2 D + \ldots + a_n D)^\star=\left(a_1 D^\star + a_2 D^\star
+ \ldots + a_n D^\star \right)^\star =J_0^\star=I^\star$, and so $I$ is
$\star$--finite.
\end{proof}

 \begin{prop} Let $\star$ be a semistar operation on an integral
domain $D$  and let $I \in \boldsymbol{\overline{F}}(D)$.  Then
$I$ is
quasi--$\star_{\!_f}$--invertible  if and only if  $I$ and $(D^{\star}:I)$
are $\star_{\!_f}$--finite
(hence, $\star$--finite)  and  $I$ is
quasi--$\star$--invertible.
\end{prop}
\begin{proof} Let $\iota:D \hookrightarrow D^\star$ be the canonical embedding
and let $\star_\iota$ be the (semi)star operation on $D^\star$ canonically
induced by $\star$.

For the ``if'' part, use the same argument of the proof of the
``if'' part of Proposition \ref{cor:invfin}.

  The  ``only if'' part.  Since $I$ is
  quasi--$\star_{\!_f}$--invertible, then $(D^{\star}:I)$ is also
  quasi--$\star_{\!_f}$--invertible, thus we have that $I^{\star_{\!_f}}$
  and $\left(D^{\star}:I\right)^{\star_{\!_f}}=(D^{\star}:I)$ are
  $(\star_{\!_f})_\iota$--invertibles.  Then, $I^{\star_{\!_f}}$ and
  $(D^{\star}:I)$ are $(\star_{\!_f})_\iota$--finite (Corollary
  \ref{cor:invfin}) and then $I$ and $(D^{\star}:I)$ are
  $\star_{\!_f}$--finite, by Lemma \ref{lemma:starifinite}.  Clearly $I$
  is quasi--$\star$--invertible, since $\star_{\!_f} \leq \star$ (Lemma
  \ref{lemma:qinv} (1)).  \end{proof}

  It is natural to ask  under which conditions a quasi--$\star$--invertible
fractional ideal is $\star$--invertible. Let $I  \in
\boldsymbol{{F}}(D)$ be
quasi--$\star$--invertible.  Then $\left(I(D^\star:I)\right)^\star=D^\star$.
Suppose that $I$ is also $\star$--invertible, that is,
$(I(D:I))^\star=D^\star$.  Then,
$(D:I)^\star=((D:I)\left(I(D^\star:I))^\star\right)^\star=(((D:I)I)^\star
\left(D^\star:I)\right)^\star=\left(D^\star:I\right)^\star =
(D^\star:I) = (D^\star : I^\star) \supseteq (D:I)^\star$.
Therefore we have the following  (cf. also Remark \ref{rem:2.2}(b)):

\begin{prop}\label{pr:2.16} Let $\star$ be a semistar operation on an integral
domain $D$.  Let $I$ be a quasi--$\star$--invertible fractional ideal of
$D$.  Then, $I$ is $\star$--invertible if and only if
$(D:I)^\star=(D^\star:I)$\ (i.e. $ \left(I^{-1}\right)^\star =
\left(I^\star\right)^{-1}$). \hfill $\Box$
\end{prop}

The following corollary is straightforward (in particular, part (2) follows
immediately from \cite[proof of Remark 1.7]{FH2000}):

\begin{cor} \label{cor:2.16}
Let $\star$ be a semistar operation on an integral domain $D$, and
let  $I \in \boldsymbol{F}(D)$.
\begin{enumerate}
\item[(1)] If $\star$ is a (semi)star operation then $I$ is
quasi--$\star$--invertible if and only if $I$ is $\star$--invertible.
\item[(2)] If $\star$ is stable and $I \in \boldsymbol{f}(D)$ then $I$
is quasi--$\star$--invertible if and only if $I$ is $\star$--invertible.
\hfill $\Box$
\end{enumerate}
\end{cor}

We notice that if $\star$ is a semistar operation of finite type,
 $\star$--invertibility depends only on the set of
quasi--$\star$--maximal ideals of $D$.  Indeed, it is clear that $I \in
\boldsymbol{F}(D)$ is $\star$--invertible if and only if
$\left(II^{-1}\right)^\star \cap D$ is not contained in any
quasi--$\star$--maximal ideal.  Then, from Lemma \ref{lemma:max}, we
deduce immediately the following  general result (cf.  \cite[Proposition
4.25]{FH2000}):

\begin{prop} \label{prop:starinvtilinv}
Let $\star$ be a semistar operation on an integral domain $D$. Let
$I \in \boldsymbol{F}(D)$. Then $I$ is $\star_{\!_f}$--invertible if and
only if $I$ is ${\tilde{\star}}$--invertible. \hfill $\Box$
\end{prop}

  A classical example due to Heinzer can be used for describing the
content of the previous proposition.

\begin{ex} \label{ex:2.19}
    Let $K$ be a field and $X$ an indeterminate over $K$.  Set $D:= K[\![X^{3},
    X^{4}, X^{5}]\!]$ and $M:=(X^{3},
    X^{4}, X^{5})D$.  It is easy to see that $D$ is a one-dimensional
    Noetherian local integral domain, with maximal ideal $M$.  Let $\star :=
    v$, note that in this case $ v = \star = \star_{_{\!f}}=t$ and
    $\mathcal{M}(v) =\{M\}$, since $ M=(D:K[\![X]\!])$.  Therefore, $w
    =\widetilde{v} = d $.  In this situation $\Inv(D, v) =\Inv(D, t)
    =\Inv(D, w)$ $=\Inv(D) =\{zD \mid\ z\in K\,,\; z\neq 0\}$.  But $v=t
    \neq w =d$, because in general $(I\cap J)^t$ is different from $I^t\cap
    J^t$ in $D$, since $D$ is not a Gorenstein domain \cite[Theorem 5,
    Corollary 5.1]{A88}  and \cite[Theorem 222]{Kaplansky70}.  \end{ex}

  A result ``analogous'' to Proposition \ref{prop:starinvtilinv}
does not hold, in general, for quasi-semistar-invertibility, as we show
in the following:

\begin{ex}\label{ex:2.20}
Let $D$ be a pseudo--valuation domain, with maximal ideal $M$, such that
$V := M^{-1}$ is a DVR (for instance, take two fields $k \subsetneq
K$ and let $V:=K[\![X]\!]$, $M := XK[\![X]\!]$ and $D :=  k+M$).  Consider
the semistar operation of finite type $\star := \star_{\{V\}}$, defined
by $E^{\star_{\{V\}}}:= EV$, for each $E\in \boldsymbol{\overline{F}}(D)$.
It is clear that $M$ is the only quasi--$\star$--maximal ideal of
$D$. Thus, ${\tilde{\star}}=\star_{\{M\}}= d$, the identity
(semi)star  operation of $D$.  We have $(M(V:M))^\star=(M(V:M))V=V$,
since $V$ is a DVR. Hence, $M$ is quasi--$\star$--invertible.  On the
other side, $M$ is not invertible (i.e., not
quasi--${\tilde{\star}}$--invertible), since $MM^{-1} =MV=M$.
\end{ex}

  Under the assumption $D^\star=D^{\tilde{\star}}$ we have the
following extension of Proposition \ref{prop:starinvtilinv} to the case of
quasi--semistar--invertibility:

\begin{prop}
Let $\star$ be a semistar operation on an integral domain $D$.
Suppose that $D^\star = D^{\tilde{\star}}$. Let $I \in
\boldsymbol{\overline{F}}(D)$.  Then $I$ is
quasi--$\star_{\!_f}$--invertible if and only if $I$ is
quasi--${\tilde{\star}}$--invertible
\end{prop}
 \begin{proof}
If $I$ is quasi--${\tilde{\star}}$--invertible, then there exists $J \in
\boldsymbol{\overline{F}}(D)$ with
$(IJ)^{\tilde{\star}}=D^{\tilde{\star}}$.  This implies
$(IJ)^{\star_{\!_f}}=D^{\star_{\!_f}}$, since ${\tilde{\star}} \leq
\star_{\!_f}$.  Conversely, suppose that there exists $J \in
\boldsymbol{\overline{F}}(D)$ such that
$(IJ)^{\star_{\!_f}}=D^{\star_{\!_f}}$.  Then $IJ \subseteq
D^{\star_{\!_f}}=D^\star=D^{\tilde{\star}}$.  Thus,
$(IJ)^{\tilde{\star}} \subseteq D^{\tilde{\star}}$.  If
$(IJ)^{\tilde{\star}} \subsetneq D^{\tilde{\star}}$, then
$(IJ)^{{\tilde{\star}}} \cap D \subsetneq D$ is a
quasi--${\tilde{\star}}$--ideal of $D$.  It follows that
$(IJ)^{\tilde{\star}} \cap D$ is contained in a
quasi--${\tilde{\star}}$--maximal $P$ of $D$.  From Lemma \ref{lemma:max},
$P$ is also a quasi--$\star_{\!_f}$--maximal.  Then, $(IJ)^{\star_{\!_f}}
\cap D \subseteq ((IJ)^{{\tilde{\star}}} \cap D)^{\star_{\!_f}}
\subseteq P^{\star_{\!_f}} \subsetneq D^{\star_{\!_f}}$, a
contradiction.  Then, $I$ is quasi--${\tilde{\star}}$--invertible.
\end{proof}

 \begin{rem} \label{rk:2.22} \bf (a) \rm If $\star$ is a semistar
operation on an integral domain $D$, we already observed (Remark
\ref{rem:2.2}(a)) that $\Inv(D, \star)$ is not a group with respect to
the standard multiplication of fractional ideals.  In the set of \sl the
$\star$--invertible $\star$--fractional ideals,\rm \ i.e. in the set
$\Inv^\star(D) := \{I \in \Inv(D, \star) \mid I = I^\star\}$,\ we can
introduce \it a semistar composition \rm ``$\times$'' in the following
way $ I\times J := (IJ)^\star$.  Note that $(\Inv^\star(D), \times)$ is
still not a group, in general, because for instance it does not possede
an identity element (e.g. when $D^\star \in \boldsymbol{\overline{F}}(D)
\smallsetminus \boldsymbol{F}(D)$).

On the other hand, $\QInv^\star(D) := \{I \in
\QInv(D, \star) \mid I = I^\star\}$, with the semistar composition
``$\times$'' introduced above, is always a group, having as identity $D^\star$
and unique inverse of $I\in \QInv^\star(D)$ the $D$--module $(D^\star
: I) \in \boldsymbol{\overline{F}}(D)$, which belongs to $\QInv^\star(D)$.
This fact provides also one of the motivations for considering $\QInv(D,
\star) $ and $\QInv^\star(D)$ (and not only $\Inv(D,
\star) $ and $\Inv^\star(D)$, as in the ``classical'' star case).

It is not difficult to prove that: \it let $\star$ be a semistar
operation on an integral domain $D$, then:
$$
(\Inv^\star(D), \times) \mbox{ is a group} \;\; \Leftrightarrow \;\;  (D :
D^\star) \neq (0)\,.
$$
 \rm

As a matter of fact, $(\Rightarrow)$ holds because $D^\star \in
\Inv^\star(D) \subseteq \boldsymbol{F}(D)$ and so $(D : D^\star) \neq
(0)$.  $(\Leftarrow)$ holds because $(D : D^\star) \neq (0)$ implies that
$D^\star \in \Inv^\star(D)$ and, for each $I\in \Inv^\star(D)$, we have
$(D^\star :I) \in \boldsymbol{{F}}(D)$, thus $(D:I)^\star =(D^\star :I)$
(Remark \ref{rem:2.12}(d1)) and so the inverse of each element $I\in
\Inv^\star(D)$ exists and is uniquely determined in $\Inv^\star(D)$.

 Note that, even if $(\Inv^\star(D), \times)$ is a group, $\Inv^\star(D)$
could be a proper subset of $\QInv^\star(D)$.  For this purpose, take $D,\ V,\
M$ as in Example \ref{ex:2.20}, in this case $D^\star = V$ and $(D:V) =
M \neq (0)$, hence $(\Inv^\star(D), \times)$ is a group, but $M \in
\QInv^\star(D) \smallsetminus \Inv^\star(D)$.

\bf (b) \rm Note that,  if $\star$ is a semistar
operation on an integral domain $D$, the group $(\QInv^\star(D),
\times)$ can be identified with a more classic group of star-invertible
star-ideals.  As a matter of fact, it is easy to see that:

\centerline{$ (\QInv^\star(D), \times) \ = \
(\Inv^{\star_{\iota}}(D^\star), \times')$}

\noindent where $\iota:D \rightarrow D^\star$ is the canonical
embedding, $\star_{\iota}$ is the (semi)star operation on $D^\star$,
canonically associated to $\star$ (Proposition \ref{stari}), and the
(semi)star composition ``$\times'$'' in $\Inv^{\star_{\iota}}(D^\star)$
is defined by $E\times' F := (EF)^{\star_{\iota}}$.

\bf (c) \rm Let $\star_{1},\ \star_{2}$ be two semistar operations on an
integral domain $D$.  If $\star_{1} \leq \star_{2}$ then $\Inv(D,
\star_{1}) \subseteq \Inv(D, \star_{2})$ and $\QInv(D,
\star_{1}) \subseteq \QInv(D, \star_{2})$.  Note that it is not true in
general that $\Inv^{\star_{1}}(D) \subseteq \Inv^{\star_{2}}(D)$ or that
$\QInv^{\star_{1}}(D) \subseteq \QInv^{\star_{2}}(D)$, because there is
no reason for a $\star_{1}$--ideal (or --module) to be a
$\star_{2}$--ideal (or --module).\   For instance, let $T$ be a proper overring
of an integral domain $D$, let $\star_{1}:= d$ be the identity (semi)star
operation on $D$ and let $\star_{1} :=\star_{\{T\}}$ be the semistar
operation on $D$ defined by $E^{\star_{\{T\}}}:= ET$, for each $E\in
\boldsymbol{\overline{F}}(D)$.  If $I$ is a nonzero principal ideal of
$D$, then obviously $I\in \Inv^{\star_{1}}(D) \ (= \Inv(D) =
\QInv^{\star_{1}}(D))$ but $I$ does not belong to $\QInv^{\star_{2}}(D)
$ (and, in particular, it does not belong to $\Inv^{\star_{2}}(D))$,
because $I^{\star_{2}}= IT \neq I$.

Note that, even if $\Inv(D, {\star_{1}}) = \Inv(D, {\star_{2}})$, for some pair of
semistar operations $\star_{1} \leq \star_{2}$), it is not true in
general that $\Inv^{\star_{1}}(D) \subseteq \Inv^{\star_{2}}(D)$.  Take
$D,\ V,\ M$ as in Example \ref{ex:2.20}.  Let $\star_{1}:= d$ be the
identity (semi)star operations on $D$ and let $\star_{2}
:=\star_{\{V\}}$.  In this case, $\Inv(D, {\star_{1}}) = \Inv(D, {\star_{2}})$,
because $\star_{1}= \widetilde{\ \star_{2}}$ and $\star_{2}=
(\star_{2})_{_{\!f}}$ (Proposition \ref{prop:starinvtilinv}).  But,
$\Inv^{\star_{2}}(D) \subsetneq \Inv^{\star_{1}}(D) =\Inv(D)$, because
$\Inv^{\star_{2}}(D) \subseteq \Inv^{\star_{1}}(D)=\Inv(D)$ since each
$\star_{2}$--ideal is obviously a $\star_{1}$--ideal, and moreover the
proper inclusion holds because, as above, a nonzero principal ideal of
$D$ belongs to $\Inv(D)$ but not to $\Inv^{\star_{2}}(D)$.

On the other hand, if $\star_{1} \leq \star_{2}$ are two \sl star
operations \rm on $D$, then it is known that $\Inv^{\star_{1}}(D) \subseteq
\Inv^{\star_{2}}(D)$, essentially because, in this case,  $I \in
\Inv^{\star_{1}}(D)$ implies that $I = I^{\star_{1}} = I^v$ and so $I =
I^{\star_{2}}$ \cite[Proposition 3.3]{A'88}.

\bf (d) \rm  Let $\star$ be a semistar operation on an integral
domain $D$, let $v(D^{\star})$ be the semistar operation on $D$ defined in
Lemma \ref{lemma:qinv}(4) and let $I, J\in \boldsymbol{F}(D)$ [respectively, $I, J\in
\boldsymbol{\overline{F}}(D)$].  Assume that $I$ is a $\star$--invertible
[respectively, quasi--$\star$--invertible] $\star$--ideal of $D$, then:
$$
\left(IJ^v\right)^\star = (I: (D:J))\,\; {\mbox{[respectively,
$\left(IJ^{v(D^{\star})}\right)^\star = (I: (D^\star:J))$]}}.$$
\hskip 12pt Recall that, since $I =I^\star$, then $(I: (D:J)) =(I:
(D:J))^\star$.  It is obvious that $IJ^v \subseteq (I:(D:J^v))= (I:
(D:J))$ and thus $\left(IJ^v\right)^\star \subseteq (I: (D:J))$.
Conversely, if $z \in (I: (D:J))$ then $z(D:J) \subseteq I$ and so
$z(D:I) \subseteq J^v$.  Therefore $z \in zD^\star =z((D:I)I)^\star\subseteq
\left(IJ^v\right)^\star$.

For the quasi--$\star$--invertible case, if $I =I^\star$, then $(I:
(D^\star:J)) =(I: (D^\star:J))^\star$ and $I =ID^\star$.  It is obvious
that $IJ^{v(D^{\star})} \subseteq (I:(D^\star:J^{v(D^{\star})}))= (I:
(D^\star:J))$ and thus $\left(IJ^{v(D^{\star})}\right)^\star \subseteq (I:
(D^\star:J))$.  Conversely, if $z \in (I: (D^\star:J))$ then
$z(D^\star:J) \subseteq I$ and so $z(D^\star:I) \subseteq J^{v(D^{\star})}$.
Therefore $z \in zD^\star =z((D^\star:I)I)^\star\subseteq
\left(IJ^{v(D^{\star})}\right)^\star$.

 \end{rem}
 \smallskip

    In the next theorem, we investigate the behaviour of a
    $\star$--invertible ideal (when $\star$ is a semistar operation) with
    respect to the localizations at quasi--$\star$--maximal ideals and in
    the passage to semistar Nagata ring.  More precisely, in the spirit of
    Kaplansky's theorem on ($d$--)invertibility \cite[Theorem
    62]{Kaplansky70}, we extend a characterization of $t$--invertibility
    proved in \cite[Corollary 3.2]{MMZ88} and two Kang's results proved in
    the star setting \cite[Theorem 2.4 and Proposition 2.6]{K89}.

\begin{thm} \label{thm:2.5}
Let $\star$ be a semistar operation on an integral domain $D$.
Assume that $\star = \star_{\!_f}$.  Let $I \in \boldsymbol{f}(D) $,
then the following are equivalent:
\begin{enumerate}
\item[(i)] $I$ is $\star$--invertible.

\item[(ii)]  $ID_Q\in\Inv(D_Q)$, for each $Q\in\mathcal{M}(\star)$ (and then $ID_Q$ is principal in
$D_Q$).

\item[(iii)]  $I\Na(D,\star)\in\Inv(\Na(D,\star))$.
\end{enumerate}
\end{thm}
\begin{proof} (i) $\Rightarrow$ (ii).  If $(II^{-1})^\star = D^\star$,
then $II^{-1} \not\subseteq Q$, for each $Q \in
\mathcal{M}(\star)$. Since $I \in \boldsymbol{f}(D) $, by flatness
we have:
$$
I^{-1}D_{Q}= (D:I)D_{Q}= (D_{Q}: ID_{Q})= (ID_{Q})^{-1}\,.
$$
Therefore, for each $Q \in \mathcal{M}(\star)$, since $II^{-1}
\not\subseteq Q$, we have:
$$
D_{Q}= (II^{-1})D_{Q}= ID_{Q}I^{-1}D_{Q} =
ID_{Q}(ID_{Q})^{-1}\,.$$

(ii) $\Rightarrow$ (iii).  From the assumption and from the proof of
(i) $\Rightarrow$ (ii), we have that $II^{-1} \not\subseteq Q$, for
each $Q \in \mathcal{M}(\star)$. Since $I \in \boldsymbol{f}(D) $,
by the flatness of the canonical homomorphism $D \rightarrow
D[X]_{N(\star)}= \Na(D, \star)$, we have:
$$
(I[X]_{N(\star)})^{-1}= (D[X]_{N(\star)} :I[X]_{N(\star)})=
(D:I)[X]_{N(\star)}= I^{-1}[X]_{N(\star)}\,.$$ Since $II^{-1}
\not\subseteq Q$, then  $(II^{-1})[X]_{N(\star)} \not\subseteq
Q[X]_{N(\star)}$, for each $Q \in \mathcal{M}(\star)$.  From
\cite[Proposition 3.1(3)]{FL03}, we deduce that:
 $$
 D[X]_{N(\star)} =(II^{-1})[X]_{N(\star)} =
 I[X]_{N(\star)}(I[X]_{N(\star)})^{-1}\,,$$
 where $I\Na(D, \star) = I[X]_{N(\star)}$.

 (iii) $\Rightarrow$ (i).  From the assumption and from the previous
 considerations, we have:
$$
D[X]_{N(\star)} =
 I[X]_{N(\star)}(I[X]_{N(\star)})^{-1}=(II^{-1})[X]_{N(\star)}\,,$$
 and thus $(II^{-1})[X]_{N(\star)} \not\subseteq
Q[X]_{N(\star)}$, for each $Q \in \mathcal{M}(\star)$.  This fact
implies that $II^{-1} \not\subseteq Q$, for each $Q \in
\mathcal{M}(\star)$. From \cite[Lemma 2.4 (1)]{FL03}, we deduce
immediately that $(II^{-1})^\star = D^\star$.
\end{proof}

\begin{cor}  \label{cor:2.22}
Let $\star$ be a stable semistar operation of finite
type on $D$, and let $I \in \boldsymbol{f}(D)$. Then, the
conditions  {\rm (i)--(iii)} of Theorem \ref{thm:2.5} are equivalent to:
\begin{enumerate}
\item[(iv)] $I$ is quasi--$\star$--invertible.
\end{enumerate}
\end{cor}
\begin{proof}
Apply Corollary \ref{cor:2.16}.
\end{proof}

\begin{rem}
It is known  \cite[Proposition 2.6]{K89} (cf.  also \cite[Section
1]{AZ93} and \cite[Section 1]{CP03}) that, \it if $\star$ is a star
operation of finite type on an integral domain $D$, an ideal $I$ of $D$
is $\star$--invertible if and only if it is $\star$--finite and locally
principal (when localized at the $\star$--maximal ideals).  \rm As a
matter of fact, by Corollary \ref{cor:invfin}, we have that, if $I$ is
$\star$--invertible, then $I$ is $\star$--finite.  Moreover,
$(II^{-1})^\star = D$ implies $II^{-1} \not \subseteq Q$, for each
$\star$--maximal ideal $Q$ of $D$.  It follows that $ID_QI^{-1}D_Q=D_Q$.
Thus, $ID_Q$ is invertible (hence, principal) in $D_Q$.  Conversely,
assume that $I^\star=J^\star$, with $J \in \boldsymbol{f}(D)$, $J
\subseteq I$.  It is clear that $I^{-1} = J^{-1}$, since
$I^v=(I^\star)^v=(J^\star)^v=J^v$, being $\star \leq v$ \cite[Theorem
34.1(4)]{Gil92}. Suppose that $I$ is not $\star$--invertible, that
is, $(II^{-1})^\star \subsetneq D$.  Then, there exists a
$\star$--maximal ideal $Q$ of $D$, such that $II^{-1} \subseteq Q$.  It
follows $QD_Q \supseteq ID_QI^{-1}D_Q = ID_QJ^{-1}D_Q=ID_Q(JD_Q)^{-1}
\supseteq ID_Q(ID_Q)^{-1}$, a contradiction, since $ID_Q$ is principal.

    We will see in a moment that the ``if'' part of a similar result for
  semistar operations does not hold, even if $I = I^\star$. More precisely,
  we can extend partially \cite[Proposition 1.1]{Gabelli89} in the
  following way:

  \it Let $I \in \boldsymbol{F}(D)$ and let $\star $ be a semistar
operation on $D$, the following properties are equivalent:
\begin{enumerate}
\item[(i)] $I$ is $\star_{_{\!  f}}$--invertible;
\item[(ii)] $(Q:I) \subsetneq (D:I)$, for each $Q \in {\mathcal M}(\star_{_{\!  f}})$;
\item[(iii)] $(Q:I) \subsetneq (D:I)$, for each $Q \in {\mathcal
M}(\star_{_{\!  f}})$ and $Q\supseteq I(D:I)$.

\noindent \hskip -24 pt  Moreover, each of the previous properties implies the
following:
\item[(iv)]
$I$ is $\star_{_{\!  f}}$--finite and $ID_Q\in\Inv(D_Q)$, for each
$Q\in\mathcal{M}(\star_{_{\!  f}})$ (and so $ID_Q$ is principal in $D_Q$).
\end{enumerate}\rm

As a matter of fact, (i) $\Rightarrow$ (ii) because  $D^\star =
\left(I(D:I)\right)^\star $  and if $(Q:I) = (D:I)$, for some $Q \in
{\mathcal M}(\star_{_{\!  f}})$, then $I(D:I) = I(Q:I)\subseteq Q $, thus
$\left(I(D:I)\right)^{\star_{_{\!  f}}}\subseteq Q^{\star_{_{\!  f}}}
\subsetneq D^\star$, hence we reach a contradiction.  (ii) $\Rightarrow$
(iii) is trivial.  (iii) $\Rightarrow$ (i): if not, $ I(D:I)\subseteq
Q$, for some $Q \in {\mathcal M}(\star_{_{\!  f}})$, thus $(D:I)
\subseteq (Q:I)$ and hence $(D:I)
= (Q:I)$, which contradicts (iii).

Finally (ii) $\Rightarrow$ (iv), because of Proposition \ref{cor:invfin}
and because for $z_{Q}\in (D:I) \smallsetminus (Q:I)$, we have $z_{Q}I \subset
D\smallsetminus Q$, and so $z_{Q}ID_{Q}= D_{Q}$, i.e. $ID_{Q}
=(z_{Q})^{-1}D_{Q}$, for each $Q \in {\mathcal M}(\star_{_{\!  f}})$.

But note that, in the semistar setting, (iv) $\not\Rightarrow$ (i), even in
case $I$ is a $\star_{_{\!  f}}$--ideal, $\star_{_{\!  f}}$--finite,
as the following example will show.  However, we can re-establish a
cha\-rac\-te\-rization in the quasi--$\star$--invertibility setting in the
following way:\ \it if $\star$ is a semistar operation of finite type
on  an integral domain $D$ and if $I \in
\boldsymbol{\overline{F}}(D)$, then  $I \in   \QInv(D, \star)  $ if and
only if $I^\star$ is $\star$--finite and $I^\star {D^\star}_M$ is
principal, for each $\star_\iota$--maximal ideal $M$ of $D^\star$.
\rm
\end{rem}

\begin{ex} Let $D$ be a valuation domain, $P$ a nonzero
nonmaximal noninvertible prime ideal of $D$ such that $D_{P}$ is a
discrete valuation domain. (For instance, if $K$ is a field and $X,\ Y$
are two indeterminates over $K$, let $D: = K+XK[X]_{(X)}+
YK(X)[Y]_{(Y)}$ and $P:= YK(X)[Y]_{(Y)}$; in this case $D$ is a
two-dimensional valuation domain, $D_{P}=K(X)[Y]_{(Y)}$ and $P=PD_{P}=
YD_{P}\supsetneq YD$.)  Set $\star:= \star_{\{P\}}$.  In this situation,
$\star =\star_{\!_f}$ and $\mathcal{M}(\star) =\{P\}$, thus $\star
=\widetilde{\star}$, i.e. $\star$ is a stable semistar operation
of finite type on $D$.  Note that $P$ is in fact a $\star$--ideal
of $D$, since $P^\star = PD_{P}=P$.  Moreover, $P^{\star}= PD_P
=tD_{P}=(tD)^\star$ for some nonzero $t \in D_{P}$, i.e. $P$ is a non
zero principal ideal in $D^{\star}=D_{P}$, since $D_{P}$ is a DVR, by
assumption.  Thus, $P$ is a $\star$--ideal, $\star$--finite and locally
principal, when localized at the quasi--$\star$--maximal ideal(s) of
$D$.  But $P$ is not $\star$--invertible , since in this situation
$(D:P) = (P:P) =D_{P}$ \cite[Proposition 3.1.5]{FHP97} and hence
$(P(D:P))^\star = (P(P:P))^\star =(PD_{P})^\star = P^\star =P$.  Note
also that, in this situation, $P$ is quasi--$\star$--invertible (because
$(P(D^\star :P))^\star = (tD_{P}t^{-1}D_{P})^\star = D_{P}=D^\star$) and
$D^\star =D_{P}= (PD_{P}:PD_{P}) = (P:P)D_{P} = (P:P)^\star$.  \end{ex}

  Next two results generalize to the semistar setting \cite[Theorem 2.12 and
Theorem 2.14]{K89}.

\begin{cor}
Let $\star$ be a semistar operation on an integral domain $D$.
Assume that $\star= \star_{\!_f}$. Let $h \in D[X], h \neq 0$, then:
$$\boldsymbol{c}(h) \in   \Inv(D, \star)  \;\; \;\;\Leftrightarrow \;\; \;\; h \Na (D,\star)=\boldsymbol{c}(h)
\Na(D,\star).$$
  In particular, $\boldsymbol{c}(h) \in \Inv(D,\star)$ if and only if
  $\boldsymbol{c}(h) \in  \QInv(D, \star)  $.
\end{cor}
\begin{proof}
The proof  of the first part of the statement is based on the
following result by D.D. Anderson \cite[Theorem 1]{A77}: If $R$ is
a ring and $h \in R[X], h \neq 0$, then $hR(X) \subseteq c(h)
R(X)$ and, moreover, the following are equivalent:
\begin{enumerate}
\item[(1)] $hR(X)=\boldsymbol{c}(h)R(X)$.
\item[(2)] $\boldsymbol{c}(h)$ is locally principal (in $R$).
\item[(3)] $\boldsymbol{c}(h)R(X)$ is principal (in $R(X)$).
\end{enumerate}

($\Rightarrow$) By Theorem \ref{thm:2.5} ((i) $\Rightarrow$ (ii)),
we have that $\boldsymbol{c}(h)D_Q$ is principal, for each $Q \in
\mathcal{M}(\star)$. Hence, $$\boldsymbol{c}(h) D_Q [X]_{N(\star)}
= \boldsymbol{c}(h)(D[X]_{N(\star)})_{QD[X]_{N(\star)}} =
\boldsymbol{c}(h)D_Q(X)$$ is principal, for each $Q \in
\mathcal{M}(\star)$. By applying Anderson's result to the local
ring $R=D_Q$, we deduce that $hD_Q(X)=\boldsymbol{c}(h)D_Q(X)$,
for each $Q \in \mathcal{M}(\star)$. The conclusion
follows from Proposition \ref{prop:nagata}, (2) and (3) \\
($\Leftarrow$) If $h \Na (D,\star)=\boldsymbol{c}(h)\Na(D,\star)$,
then by localization we obtain that
$hD_Q(X)=\boldsymbol{c}(h)D_Q(X)$, for each $Q \in
\mathcal{M}(\star)$ (Proposition \ref{prop:nagata} and
\cite[Corollary 5.3]{Gil92}). By Anderson's result, we deduce that
$\boldsymbol{c}(h)D_Q$ is principal, i.e. $\boldsymbol{c}(h)D_Q
\in \Inv(D_Q)$, for each $Q \in \mathcal{M}(\star)$. The
conclusion follows from Theorem \ref{thm:2.5} ((ii) $\Rightarrow$
(i)).

  The last part of the statement follows from the fact that
$\Na(D,\star) =\Na(D,\widetilde{\star})$ \cite[Corollary 3.5(3)]{FL03}
and from Corollary \ref{cor:2.16} and Proposition
\ref{prop:starinvtilinv}\ or, directly, from Corollary \ref{cor:2.22}.
\end{proof}

\begin{prop}
Let $\star$ be a semistar operation on an integral domain $D$. If
$H$ is an invertible ideal of $\Na(D,\star)$, then $H$ is
principal in $\Na(D,\star)$.
\end{prop}
\vskip -4pt \begin{proof} We can assume that $H \in \Inv(\Na(D,\star))$ and $H \subseteq \Na(D,\star)$, then, in
particular, $H=(h_1,h_2, \ldots, h_n) \Na(D,\star)$, with $h_i \in
D[X],\ 1 \leq i \leq n$. For each $Q \in \mathcal{M}(\star_{\!_f})$, by
localization, we obtain that $HD_Q(X)=(h_1,h_2,\ldots, h_n)D_Q(X)$
is a nonzero principal ideal (Theorem \ref{thm:2.5} ((iii)
$\Rightarrow$ (ii)). By a standard argument, if $d_i:= \deg(h_i)$,
for $1 \leq i \leq n$, and if
$$h:=h_1+h_2 X^{d_1 +1}+h_3X^{d_1+d_2+2}+ \ldots + h_n X^{d_1 +
d_2 + \ldots + d_{n-1}+n-1} \in D[X], $$ then it is not difficult
to see that $HD_Q(X)=hD_Q(X)$, for each $Q \in
\mathcal{M}(\star_{\!_f})$. From Proposition \ref{prop:nagata}(3), we
deduce that $H \Na (D,\star)= h \Na(D,\star)$.
\end{proof}




\bigskip




\begin{thebibliography}{10}

\bibitem{A77}
D.~D. Anderson, \emph{Some remarks on the ring {$R(X)$}}, Comment.
Math. Univ.
  St. Paul. \textbf{26} (1977/78), no.~2, 137--140.

  \bibitem{A88} D.~D. Anderson, \emph{Star-operations induced by
 overrings},\; Comm.  Algebra \bf 16 \rm (1988), 2535--2553.


\bibitem{AC2000}
D.~D. Anderson and S.~J. Cook, \emph{Two star--operations and their
induced lattices}, Comm.  Algebra \textbf{28} (2000), no.~5, 2461--2475.


   \bibitem{AMZ89} D.~D. Anderson , Joe Mott and Muhammad Zafrullah,
  \emph{Some quotient based statements in multiplicative ideal theory},
  Boll.  Unione Mat.  Ital.  \textbf{3-B} (1989), 455--476.

\bibitem{AZ93}
D.~D. Anderson and Muhammad Zafrullah, \emph{On
{$t$}--invertibility. {III}},
  Comm. Algebra \textbf{21} (1993), no.~4, 1189--1201.


  \bibitem{A'88} David F. Anderson, \emph{A general theory of class group},\;
  Comm.  Algebra \bf 16 \rm (1988), 805--847.


 \bibitem{A00}David F. Anderson, \emph{The class group and the local
 class group of an integral domain}, ``Non-Noetherian Commutative Ring Theory'' (Scott~T. Chapman and
  Sarah Glaz, eds.), Kluwer Academic Publishers, 2000, pp.~33--55.~

   \bibitem{B82} Alain Bouvier, \emph{Le groupe des classes d'un anneau int\`egre},
 107--\`eme Congr\'es National des Soci\'et\'es Savantes, Brest, Fasc. IV (1982), 85--92.


  \bibitem{B83} Alain Bouvier, \emph{The local class group of a Krull
 domain},\;  Canad.  Math.  Bull.  \bf 26 \rm (1983), 13--19.


  \bibitem{BZ88} Alain Bouvier and Muhammad Zafrullah, \emph{On some
 class groups of an integral domain}, Bull.  Soc.  Math.  Gr\`ece \bf 29
 \rm (1988), 45--59.

 \bibitem{CP03}
Gyu~Whan Chang and Jeanam Park, \emph{Star--invertible ideals of integral
domains}, Boll.  Unione Mat.  Ital.   \textbf{6-B} (2003), no.~1,
141--150.

 \bibitem{ElBF03} Said {El~Baghdadi} and Marco Fontana, \emph{Semistar linkedness and flatness,
  {P}r\"ufer semistar multiplication domains}, Comm.Algebra \textbf{32} (2004),
  no.~3, 1101--1126.


\bibitem{FH2000}
Marco Fontana and James~A. Huckaba, \emph{Localizing systems and
semistar
  operations}, ``Non-Noetherian Commutative Ring Theory'' (Scott~T. Chapman and
  Sarah Glaz, eds.), Kluwer Academic Publishers, 2000, pp.~169--198.

\bibitem{FHP97}
Marco Fontana, James~A. Huckaba, and Ira~J. Papick, \emph{Pr\"ufer
domains},
  Marcel Dekker Inc., New York, 1997.



 \bibitem{FJS03b}{Marco Fontana, Pascual Jara and Eva Santos},
\emph{Pr{\"u}fer $\star$--multiplication domains and semistar
operations}.  J. Algebra Appl.  \bf 2 \rm (2003), 21--50.

   \bibitem{FJS03}{Marco Fontana, Pascual Jara and Eva Santos},
\emph{Local--global properties for semistar operations}, Comm.  Algebra
(to appear).


  \bibitem {FL01a} M. Fontana and K. A. Loper, \it Kronecker function rings:
  a general approach,
   \ \rm \rm in  ``Ideal Theoretic Methods in Commutative Algebra" (D.D.
   Anderson and I.J. Papick, Eds.), M. Dekker Lecture Notes Pure Appl.
  Math.  \bf 220 \rm (2001), 189--205.

   \bibitem{FL01b} M. Fontana and K. A. Loper, \it A Krull-type theorem for
  the semistar integral
  closure of an integral domain, \ \rm ASJE Theme Issue ``Commutative
  Algebra'' \bf 26 \rm
  (2001), 89--95.




\bibitem{FL03}
Marco Fontana and K.~Alan Loper, \emph{Nagata rings, {K}ronecker
function rings, and related semistar
  operations}, Comm. Algebra \textbf{31} (2003), no.~10, 4775--4805..

\bibitem{FP03}
Marco Fontana and Mi~Hee Park, \emph{Star operations and
pullbacks}, J.Algebra
  \textbf{274} (2004), no.~1, 387--421.

 \bibitem{Gabelli89} Stefania Gabelli, \emph{Completely integrally
closed domains and $t$--ideals} Boll.  Unione Mat.  Ital.  \bf 3-B \rm
(1989), 327--342.

\bibitem{GH97}
Stefania Gabelli and Evan Houston, \emph{Coherentlike conditions
in pullbacks},
  Michigan Math. J. \textbf{44} (1997), no.~1, 99--123.


\bibitem{Gil90}
Robert Gilmer, \emph{Pr\"ufer domains and rings of integer-valued
polynomials},
  J. Algebra \textbf{129} (1990), no.~2, 502--517.


\bibitem{Gil92}
Robert Gilmer, \emph{Multiplicative ideal theory},  Marcel Dekker, New
York, 1972.


 \bibitem{GMZ94} Robert Gilmer, Joe Mott and Muhammad Zafrullah,
\emph{On $t$--invertibility and comparability}, ``Commutative Ring
Theory'' ( P.-J.
Cahen, D.L. Costa, M. Fontana and S.-E. Kabbaj, eds.), Lecture Notes in
Pure and Applied Mathematics, \bf 153 \rm (1994), 141--150, Marcel Dekker,
New-York.


\bibitem{GV77}
Sarah Glaz and Wolmer~V. Vasconcelos, \emph{Flat ideals. {II}},
Manuscripta
  Math. \textbf{22} (1977), no.~4, 325--341.

   \bibitem{HK98} Franz Halter-Koch, \emph{Ideal systems.  An
  introduction to multiplicative ideal theory}, Marcel Dekker, New-York, 1998.


   \bibitem{HK01} Franz Halter-Koch, \emph{Localizing systems, module
  systems and semistar operations}, J. Algebra \bf 238 \rm (2001),
  723--761.

    \bibitem{HK03} {Franz Halter-Koch}, \emph{Characterization of Pr\"ufer multiplication
 monoids and domaind by means of spectral module theory},
 Monatsh. Math. \bf 139 \rm (2003), 19--31.

    \bibitem{HH80} J.R. Hedstrom and Evan Houston, \emph{Some remarks on
   star operations}, J. Pure Appl.  Algebra \bf 18 \rm (1980), 37--44.

 \bibitem{HZ88} Evan Houston and Muhammad Zafrullah, \emph{Integral
domains in which each $t$--ideal is divisorial}, Michigan Math.  J.
\textbf{35} (1988), 291--300.

 \bibitem{HZ89} Evan Houston and Muhammad Zafrullah, \emph{On
$t$--invertibility, II}, Comm.  Algebra \textbf{17} (1989), 1955--1989.


 \bibitem{Jaffard60} Paul Jaffard, \emph{Les Syst\`emes d'Id\'eaux},
Traveaux et Recherches Math\'ematiques, Dunod, Paris, 1960.


\bibitem{K89}
B.~G. Kang, \emph{Pr\"ufer {$v$}--multiplication domains and the
ring {$R[X]\sb
  {N\sb v}$}}, J. Algebra \textbf{123} (1989), no.~1, 151--170.

   \bibitem{Kaplansky70} Irving Kaplansky, \emph{Commutative rings},
  Allyn and Bacon, Inc., Boston, 1970.

  \bibitem{MMZ88}
 S. Malik, Joe Mott and Muhammad Zafrullah, \emph{On $t$--invertibility},
 Comm.  Algebra \bf 16 \rm (1988), 149--170.

  \bibitem{M98} Ry{\=u}ki  Matsuda, {\it Kronecker function rings of semistar
 operations on rings},
  Algebra Colloquium {\bf 5} (1998), 241--254.

   \bibitem{M00} Ry{\=u}ki Matsuda, \emph{Note on vluation rings and
  semistar operations}, Comm.  Algebra \bf 28 \rm (2000), 2515--2519.

     \bibitem{MS96} Ry{\=u}ki Matsuda and I. Sato, {\it Note on star
    operations and semistar operations}, Bull.  Fac.  Sci.  Ibaraki Univ.
    Ser.  A {\bf 28} (1996), 5--22.

    \bibitem{MSu95} Ry{\=u}ki Matsuda and T. Sugatani, {\it Semistar
   operations on integral domains}, II, \ Math.  J. Toyama Univ.  {\bf 18}
   (1995), 155--161.

     \bibitem{Mi03} A. Mimouni, \emph{Semistar operations of finite
    character on integral domains}, Preprint.

      \bibitem{MiSa03} A. Mimouni and M. Samman \emph{Semistar operations
     on valuation domains}, Int.  J. Comm.  Rings \bf 2 \rm (2003), (to appear).


      \bibitem{O03} Akira Okabe, \emph{Some results on semistar
     operations}, Int.  J. Comm.  Rings (to appear).


\bibitem{OM94}
Akira Okabe and Ry{\=u}ki Matsuda, \emph{Semistar--operations on integral
domains}, Math.  J. Toyama Univ.  \textbf{17} (1994), 1--21.

  \bibitem{OM97} Akira Okabe and Ry{\=u}ki Matsuda, {\it Kronecker
 function rings of semistar ope\-ra\-tions}, Tsukuba J. Math., {\bf 21}
 (1997), 529--540.

    \bibitem{Q} Julien Querr\'e, \emph{Sur une propriet\'e des anneaux
  de Krull}, Bull.  Sci. Math.  \bf 95 \rm (1971), 341--354.

    \bibitem{WM97} Wang Fanggui and R. L. McCasland, \emph{On
   $w$--modules over strong Mori domains}, Comm.  Algebra  \bf 25 \rm (1997),
   1285--1306.

      \bibitem{WM97b} Wang Fanggui and R. L. McCasland, \emph{On
   $w$--projective modules and $w$--flat modules},  Algebra Colloquium \bf 4
   \rm (1997), 111--120.


         \bibitem{WM99} Wang Fanggui and R. L. McCasland, \emph{On
        strong Mori domains}, J. Pure Appl.  Algebra  \bf 135 \rm (1999), 155-165.



  \bibitem{Z1} Muhammad Zafrullah, \emph{Ascending chain conditions and
  star operations}, Comm. Algebra \bf 17 \rm (1989), 1523--1533.

  \bibitem{Z2} Muhammad Zafrullah, \emph{Putting $t$--invertibility to use},
  ``Non-Noetherian Commutative Ring Theory'' (Scott~T. Chapman and
  Sarah Glaz, eds.), Kluwer Academic Publishers, 2000, pp.~429--457.


\end{thebibliography}
\end{document}